\documentclass[11pt,a4paper]{article}
\usepackage[novbox]{pdfsync}
\usepackage{enumerate}
\usepackage{enumitem}
\usepackage{tabularx}
   \usepackage[margin=0.5in]{geometry}
\usepackage[utf8]{inputenc}
\usepackage{graphicx}
\graphicspath{{Figures/}}
\usepackage[T1]{fontenc}
\usepackage{fancyhdr, float,fancybox}
\usepackage{mathrsfs,amssymb,amsmath,amsfonts,booktabs,array,xcolor,url,cite,
color,multirow,multicol}
\usepackage{slashbox}
\usepackage{stfloats}
\usepackage{extarrows}
\usepackage{textcomp}
\usepackage{ulem}
\usepackage{cases}
\usepackage{bm}
\usepackage{arydshln}
\usepackage{hyperref}\usepackage{array,subcaption}\usepackage{verbatim}
\usepackage[font=footnotesize]{caption}
\usepackage[all]{xy}
\usepackage{listings}
\usepackage{tikz, pgf,graphicx, pgfplots, color, xcolor}
	\pgfplotsset{compat=1.12} 
\usepackage{tikz-3dplot}
\usepackage{tikz-cd}
	\usetikzlibrary{positioning}

\newtheorem{theorem}{Theorem}
\newtheorem{remark}{Remark}
\newtheorem{proposition}{Proposition}
\newtheorem{definition}{Definition}
\newtheorem{corollary}{Corollary}
\newtheorem{example}{Example}
\newtheorem{lemma}{Lemma}
\newtheorem{assumption}{Assumption}
\newcommand{\bq}{\boldsymbol{q}}
\newtheorem{proof}{Proof}

 \def\eeD{\end{definition}} \def\beD{\begin{definition}}
\def\beR{\begin{remark}} \def\eeR{\end{remark}}
\def\beL{\begin{lemma}} \def\eeL{\end{lemma}}
\def\beC{\begin{corollary}
  }\def\eeC{\end{corollary}}
  \def\beT{\begin{theorem}}\def\eeT{\end{theorem}}
  \def\beP{\begin{proposition}} \def\eeP{\end{proposition}}
\def\beXa{\begin{example}} \def\eeXa{\end{example}}
\def\beA{\begin{assumption}} \def\eeA{\end{assumption}}
\def\im{\item}\def\com{compartment}     \def\vi{\overset{\rightarrow}{\i}}   \def\lab{\label}
\def\D{\Delta} 

   \def\ei{a_i}  \def\eig{eigenvalue} 
    \def\NGM{next generation matrix}
\def\brn{basic reproduction number}\def\DFE{disease-free equilibrium}
    \def\rd{\r_{dfe}} \def\sd{\s_{dfe}}

    \def\va{\vec \al } 

\renewcommand{\i}{\;\mathsf i}
\renewcommand{\r}{\; \mathsf r}

\def\BEN{\begin{enumerate}}  \def\BI{\begin{itemize}}
\def\EEN{\end{enumerate}}   \def\EI{\end{itemize}}
\newcommand{\beq}{\begin{eqnarray}
    }
\def\eeq{\end{eqnarray}}
   \newcommand{\be}[1]{\begin{equation}\label{#1}}
\newcommand{\ee}{\end{equation}}
\def\bea{\begin{eqnarray*}} \def\im{\item} \def\Lra{\Longrightarrow}  \def\eqr{\eqref}  
\def\no{\nonumber} \def\mS{{\mathcal S}}
\newcommand{\mC}{\mathcal{C}}
\newcommand{\R}{\mathbb{R}}
\newcommand{\N}{\mathbb{N}}

\newcommand{\A}{\mathcal{A}}
 \def\sd{\s_{dfe}}\def\al{\al}
\def\b{\beta} \def\g{\gamma}    \def\de{\delta}\def\d{\delta}
\def\z{\zeta}   
  
   \def\al{\al}
 \def\La{\Lambda}\def\mR{\mathcal R} \def\fr{\frac}
\def\bc{\begin{cases}
  }
\def\ec{\end{cases}}
\def\bea{\begin{eqnarray*}} 
  
   \def\satg{satisfying } \def\Mp{More precisely, }   \def\saty{satisfy}
\def\eea{\end{eqnarray*}} \def\T{\widetilde}
   \def\for{\forall}\def\I{\infty}
\def\no{\nonumber}  
\long\def\symbolfootnote[#1]#2{
\begingroup
\def\thefootnote{\fnsymbol{footnote}}\footnote[#1]{#2}
\endgroup}

\providecommand{\pp}[1]{\left[#1\right]} 
\providecommand{\pr}[1]{\left(#1\right)} 
\def\bep{\begin{pmatrix}} \def\eep{\end{pmatrix}}
\def\bev{\begin{vmatrix}} \def\eev{\end{vmatrix}}
\def\Fr{Furthermore, }
\def\wrt{with respect to }
  \def\resp{respectively}

\newcommand\CRN{chemical reaction network }

\newcommand{\bb}{\mathbf{b}}

\newcommand{\xb}{\mathbf{x}}
\newcommand{\yb}{\mathbf{y}}
\newcommand{\Yb}{\mathbf{Y}}

\newcommand{\Wb}{\mathbf{W}}

\def\m0{{\mathcal R}_0} \def\la{\label}
\def\Eq{\Leftrightarrow}\def\qu{\quad}

\newcommand{\mL}{\mathcal{L}}

\newcommand\s{\boldsymbol{s}}

\renewcommand{\epsilon}{\varepsilon}	
\renewcommand{\tilde}[1]{\widetilde{#1}}

\newcommand\MM{{\mathcal M}}

\renewcommand{\cdot}{}

\renewcommand{\tilde}{\widetilde}

\newcommand*{\Scale}[2][4]{\scalebox{#1}{$#2$}}%
\def\frt{furthermore }

\DeclareMathOperator\ch{char} %

\newcommand{\inD}[1][\relax]{\def\argone{#1}\def\temprelax{\relax}
\ifx\argone\temprelax\right.\else\,\middle|#1\right.{}\fi}

\allowdisplaybreaks

\def\bb{\bff \beta} \newcommand{\bff}[1]{{\mbox{\boldmath$#1$}}}
\def\bzn{basic replacement number}  \def\ith{it holds that }
\def\s{\;\mathsf s}\def\al{\alpha}
\def\a{\;\mathsf a} 
\renewcommand{\i}{\;\mathsf i}
\renewcommand{\r}{\; \mathsf r}
\def\m{b}\def\bb{\bff \beta}\def\ME{mathematical epidemiology}
\def\AB{$(A,B)$ Arino-Brauer epidemic models}
\def\mG{\mathcal G}

\def\La{\Lambda}\def\m{\Lambda}\def\z{\;\mathsf z}\def\ch{characteristic polynomial }

\def\fp{fixed point}
\def\sm{stoichiometric matrix}\def\MTT{matrix tree theorem}
\def\MM{Michaelis--Menten}

\def\bb{\bff \beta} 
\def\wrt{with respect to }\def\ch{characteristic polynomial}
\def\FHJ{Feinberg--Horn--Jackson}
\def\Oth{On the other hand, }
\def\ACR{absolute concentration robustness}
\newtheorem{open}{Open Problem}
\def\eeO{\end{open}} \def\beO{\begin{open}}
\def\ME{mathematical epidemiology}\def\SPF{SIR-PH-FA model}
\def\NGM{next generation matrix}

\renewcommand{\r}{ \mathsf r}\def\mS{{\mathcal S}}\def\mC{{\mathcal C}}\def\mR{{\mathcal R}}\def\Ma{{Mathematica}}\def\mN{{\mathcal N}}
\def\GMAK{generalized mass-action kinetics}

\begin{document}
\title{Advancing
 Mathematical Epidemiology and Chemical
Reaction Network Theory via Synergies Between Them}
	\author{
 Florin Avram$^{1}$,  Rim Adenane$^{2}$, Mircea Neagu$^{3}$,
}
\maketitle

\begin{center}
	
		$^{1}$
		Laboratoire de Math\'{e}matiques Appliqu\'{e}es, Universit\'{e} de Pau, 64000, Pau,
 France; avramf3@gmail.com \\
$^{2}$ \quad Laboratoire des Equations aux dérivées partielles, Alg\'ebre et G\'{e}om\'{e}trie spectrales, d\'epartement des Math\'ematiques, Universit\'e Ibn-Tofail, 14000, Kenitra,
 Maroc; rim.adenane@uit.ac.ma \\
 $^{3}$ \quad  Transilvania University of Bra\c{s}ov, Department of Mathematics and Computer Science, 500091, Bra\c{s}ov, Romania; 	mircea.neagu@unitbv.ro
	\end{center}

\begin{abstract}

Our paper reviews some key concepts in chemical reaction network theory and  mathematical epidemiology, and examines their intersection, with  three goals.
The first is to make the case that  mathematical epidemiology (ME), and also related sciences like population dynamics, virology, ecology, etc., could benefit  by adopting the universal language of essentially non-negative kinetic systems as developed by  \CRN\ (CRN) researchers.
In this direction, our investigation of the relations between CRN and ME lead us to propose for the first time   a definition of ME models,  stated in Open Problem 1.
Our second goal is to inform  researchers outside ME of the convenient \NGM\ (NGM) approach for studying the stability of boundary points,  which do not seem sufficiently well known.
 Last but not least, we  want to help students and  researchers who know nothing about either ME or CRN to learn them quickly, by offering them a Mathematica package ``bootcamp'', {located at} \url{https://github.com/adhalanay/epidemiology_crns},
  including illustrating notebooks (and certain sections below will contain
 associated suggested notebooks; however, readers with experience may safely skip the bootcamp).  We hope that the files indicated in the titles of various sections will be helpful, though of course improvement  is always possible, and we ask the help of the readers for that.

 \end{abstract}


\tableofcontents

\section{{Introduction}}
\subsection{Motivation}
Dynamical systems  are  very important in the  ``sister sciences'' of mathematical epidemiology (ME), virology, ecology, population dynamics, \CRN\  (CRN) theory, and other related domains. They may have  very complex behaviors, including in two dimensions (as illustrated by Hilbert's 16th problem, for example).

It has long been noted, however, that the natural restriction to essentially non-negative (i.e., positivity preserving)
 ``mass-action kinetics'' leads often to   results which are  surprisingly simple,
 despite high-dimensionality (see for example Gunawardena \cite{Gun}).  This has motivated
 several researchers, already long ago,  to propose turning chemical reaction networks  theory (CRNT) (which grew out of mass-action polynomial kinetics theory) into a unified tool for  studying all applied disciplines involving essentially non-negative dynamical systems.

  There are not many  sources which have attempted to develop a unified view point of essentially non-negative systems. For one exception, see   the book \cite{Haddad}, and recent  papers which might be associated to  the new unifying banner of ``algebraic biology''---see, for example  \cite{pachter2005algebraic,macauley2020case,torres2021symbolic}.

 Unfortunately, the opposite of unification is also happening. The sister sciences seem to be growing further and  further apart, due  to their focus on particular examples, to the point that a suspicious reader might ask themselves whether sometimes they might not be studying the same system  under different names, without being aware of each other's results.

 One striking example of this is the fundamental equivalence  of mass-action polynomial kinetics to the ``absence of negative cross-terms'', already discovered by Vera H\'ars and J\'anos T\'oth \cite{hun}, which seems unknown outside \CRN\ theory, and it has been reproved in particular examples an uncountable number of times.

We will start now our campaign for unification by  short introductions to ME and CRN.
\subsection{A Short History of Mathematical Epidemiology}
{The role of the father of epidemiology has been assigned by \cite{Jones68} to the Greek physician Hippocrates (460--377 B.C.E), who first described the connection between disease and environment.
The first mathematical model of epidemiology for smallpox was formulated and solved  in 1766 by the Swiss
mathematician Daniel Bernoulli  \cite{Bernoulli1766}. In modern times, Ronald Ross  pioneered applying mathematical methods to  the study of malaria \cite{ross1916}.

In 1927, Kermack and Mckendrick proposed the foundational  SIR model which consists of a division of the population into compartments of Susceptible--Infected--Removed individuals \cite{Ker}. Their model assumes that the population is constant and an interaction is allowed between the whole classes of the population. The first non-polynomial SIR is due to
Bailey  \cite{Bailey75}, who replaced the linear force of infection in SIR by a fractional one,   given by the fraction represented by the  infected out of the sum of the susceptible and infected (i.e., $I/(S+I$).
After this era, a panoply of work on generalizations which became known as  compartmental models was carried out---see, for example, \cite{dietz1988,Diek,AndMay_bk, Heth, diekmann2000mathematical,Van, Brauer05, arino2006simple,Van08, BrauerWu,hethcote2009epidemiology,DiekHeesBrit}.

A big part of these works was dedicated to studying  the stability of the \DFE\ (DFE---see below), which is the most important concept of \ME, and by relating it to an important parameter known as the basic reproduction number $R_0$, whose origins come from population dynamics and the theory of branching processes. This research culminated into   the celebrated new generation matrix \mbox{approach~\cite{Diek,Van,Van08}}.

\subsection{The {Stoichiometric Matrix} for CRN Representation of  a SIRS Model: SIRSclosed.nb}

We must clarify from the start that CRNT  has been very efficiently summarized in tutorials by  ``great masters''  like   Gunawardena, Angeli, Craciun and D.A. Cox \cite{Gun,Ang,Yu,Cox}---see also \cite{horn1972,feinberg1987,Toth,feinberg2019foundations} for some original results, and for more extensive treatments.  Since we cannot  compete with their works, we urge the reader to peruse these sources, among others.

We will diverge now  from these works by introducing  CRNT via a  SIRS ODE {example}:
\be{SIRS0}
x'=\bep
s' \\ i' \\ r' \eep =  \bep - 1&0&1&-1 \\
  1&- 1&0&0\\
  0&1&- 1&1
\eep \bep
\beta s  i  \\ \g_i i   \\ \g_r r \\ \g_s s  \eep:=\Gamma {\bf R}(x)
\ee
(see below for the practical meaning of the variables and parameters).

Above, following the CRN literature, we have defined the model  by a triple $(\Gamma, {\bf R}, x)$ so that
 \be{crnsys}
 x'=\Gamma .{\bf R}(x)
 \ee
where we have the following:
  \BEN \im $x$ denotes the variables (called species);\im
   $\Gamma$ is the ``stoichiometric matrix'' (SM), whose columns  represent directions in which several species/\com s change simultaneously;
   \im  ${\bf R}(x)$ is the  vector of rates of change associated to each direction, also known as kinetics, which must satisfy the following admissibility conditions: \begin{enumerate}
 		\item[\bf {A1.}] $R_j(x)$ is continuously differentiable, $j=1,..,n_r$;
 		\item[\bf {A2.}] if $ \Gamma_{ij}< 0$, then $x_i=0$ implies  $R_j(x)=0$; %
 		\item[\bf {A3.}] ${\partial R_j}/{\partial x_i}(x) \ge 0$ if $ \Gamma_{ij}< 0$ and ${\partial R_j}/{\partial x_i}(x)\equiv 0$ if $ \Gamma_{ij}\ge 0$;
 		\item[\bf {A4.}] The inequality in   {A3} holds strictly for all positive concentrations, i.e., when $x \in \mathbb R_+^n$.
 	\end{enumerate}

    \EEN
\beR
 \lab{r:str}
The representation \eqr{crnsys} suggests that the investigation of the ODE of CRN might be  split in two parts: ``algebraic structural interaction elements'' embodied in $\Gamma$,  and ``parametric analysis features'' reflected in the specific structure of the functions $ \bf R(x)$ (polynomial, or Michaelis--Menten, or Hill  kinetics, etc.).
This opened up the fruitful investigation of ``robust questions'', which are independent of the specific kinetics---see, for example, ``robust stability'' \cite{AASComp,Ali23gr}---and was also exploited in ME recently in \cite{AAV}. \eeR

\subsection{Polynomial ODEs, Essential Non-Negativity and Mass-Action Representation}\lab{s:MAK}
When the rates ${\bf R}(x)$ are monomials, an alternative representation for  \eqr{crnsys} is
\beq \lab {Wdef} && \dot{\xb}=f(\xb)=\sum_{k=1}^{n_R} {\mathbf{w}}_r \xb^{\yb_r}=\Wb \xb^{\Yb}, \quad  \xb, \yb_r, {\mathbf{w}}_r \in \R_+^{n \times 1}, \eeq
where  $\Yb \in \R^{n \times n_R}$ is the ``matrix of  exponents'',
  $\Wb\in \R^{n \times n_R}$ is the ``matrix of  direction vectors'' (formed \resp\ by joining exponents $\yb_1,\dots ,\yb_{n_R}$ and directions ${\mathbf{w}}_1,\dots ,{\mathbf{w}}_{n_R}$ as columns), and where $\xb^{\Yb} \in \R_+^{n_R \times 1}$ is a column vector of monomials, whose $r$'th component is $\xb^{\yb_r}$, with $\xb^{\yb_r}:=x_1^{{\yb_r}_{1}} x_2^{{\yb_r}_{2}}\dots \in \R_+$.

   Note that any polynomial dynamical systems can be uniquely written in such a form for some distinct $\yb_i$ and
 non-zero ${\mathbf{w}}_i$ \cite{CraJinYu},  but $\Yb$ is not unique since its  dimension may be easily increased.

 Note also the ``pseudo-linearity'' property of \eqr{crnsys} and \eqr{Wdef} which may be transformed into linearity by selecting the reaction rates as variables.

\beR The fact that the parametrization \eqr{Wdef} employs only two matrices $(\Wb , {\Yb})$ for
describing any polynomial ODE probably explains why several researchers
have proposed in the past using the language of \CRN\ theory as a unifying modeling tool for all the ``sister disciplines'' which study ``essentially non-negative ODEs'' (see next section):   mathematical epidemiology, ecology, virology, biochemical systems, etc.
   \eeR
\subsection{Essentially Non-Negative Kinetic Systems}
 Kinetic systems (which is the old,  physics-inspired name for CRN systems) must be ``essentially non-negative'',   meaning that they leave invariant the non-negative orthant.

\beR An obvious sufficient condition for the essential non-negativity (i.e., the preservation of the
non-negative octant) of a polynomial system $X'=f(X)$ is that each component $f_i(X)$ may be decomposed
   as \begin{equation}\label{Hc} f_i(X)=g_i(X)-x_i h_i(X),\end{equation}
  where $g_i,  h_i$ are polynomials with non-negative coefficients, i.e., if all negative terms in an equation  contain the variable whose rate is  given by the equation.

   \beD Terms which do not satisfy \eqr{Hc} are called negative cross-effects. \eeD

\eeR

  \beXa The  Lorentz system (a famous example of chaotic behavior)
$$\bc x'=\sigma(y-x)\\
y'=\rho x- y-\boxed{x z}\\z'=xy -\beta z\ec$$
 does not satisfy \eqr{Hc}, due to the $-x z$ term in the $y$ equation.
 \eeXa

 The following result, sometimes called the ``Hungarian lemma'' is well known in  the \CRN\ literature  \cite{hun},  (\cite{Toth}, {Thm 6.27}):

 \beL  \la{l:hun} A polynomial system admits an essentially non-negative ``mass-action'' representation (see next section) {if and only if}
 there are no negative cross-effects, i.e., if \eqr{Hc} {holds.} 
\eeL


\subsection{The Traditional Reactions Representation of  CRN Theory, and the \FHJ\ (FHJ) Graph: SIRS.nb}

  There exists a  third parametrization used in CRN, in which each ``reaction'' (associated to a column of $\Gamma$) is represented
  as a  directed pair (input $\to$ output), in the style
    (``I'' $\to$ ''R'',\dots )

    This is formalized by introducing the following triple (to which one may associate both an ODE, and several interesting graphs; see below).

    \beD \lab{d:CRN} A  CRN is defined by a triple $\{\mathcal{S},\mathcal{C},\mathcal{R}\}$:

\begin{align*}
    \mathcal{S} &= \{S_1, \dots  , S_i,\dots , S_{|\mathcal{S}|}\}&\text{species}\\
    \mathcal{C} & = \{y_1,\dots ,y_\alpha,\dots ,y_{|\mathcal{C}|} :y_\alpha \in \mathbb{N}^{|\mathcal{S}|}\} &\text{complexes}\\
    \mathcal{R} &= \{  y_{\alpha} \xrightarrow{[k_{\{y_{\alpha} \to y_{\beta}\}} * x^{y_{\alpha,\beta}}]} y_{\beta}  : k_{\{y_{\alpha} \to y_{\beta}\}} \geq 0\}&\text{reactions},
\end{align*}
where Roman letters ($i,j$) and Greek letters ($\alpha,\beta$) are used to denote species and complex indices, respectively, and where $x=\{s_1, \dots  , s_i,\dots , s_{|\mathcal{S}|}\} $ denotes the vector of unknowns.

Finally, if the rate exponents $y_{\alpha,\beta}$ intervening in $\mR$ coincide with the input exponents $y_{\alpha}$, we have a {\bf {mass-action}} CRN (MAS/MAK), and otherwise,
~we have a {\bf {generalized mass-action}} CRN (GMAS/GMAK).

\eeD

  \beXa {  The {\bf {mass-action}} SIRS ODE \eqr{SIRS0} is induced by the reactions}

\be{RNSIR} \bc   S + I \xrightarrow{[\beta s i]} 2I\\  I \xrightarrow{[\gamma_i   i]} R  \\
 R \xrightarrow{[\g_{r} r]} S  \\
   S  \xrightarrow{[\g_{s} s]} R  \ec \ee

In this example, \begin{align*}
    \mathcal{S} &= \{S,I,R\}, x=(s,i,r),\\
    \mathcal{C} & = \{y_1=(1,1,0),y_2=(0,2,0),y_3=(1,0,0),y_4=(0,1,0),y_5=(0,0,1) \}\\
    \mathcal{R} &= \{  y_1 \xrightarrow{[\beta s i]} y_2,  y_4 \xrightarrow{[\g_i  i]} y_5,\dots \}\end{align*}

\beR   While the last three ``transfer reactions'' in \eqr{RNSIR} require no explanation, to understand why the first reaction is written as ``S+ I {$\to$} 2 I'' rather than ``S {$\to$} I''
   requires recalling that the mass-action assumption forces the coefficients appearing in the left ``reactants'' complex to coincide with the exponents of $x$ in
   the  rate $\beta s i$.
\eeR

\eeXa

   \beR  In CRN terminology, the pair, called reaction, is viewed as a transition from
    a ``source complex''/input $y_{\alpha}$ to a ``product complex''/output $y_{\beta}$.

    From the representation \eqr{RNSIR},  one may construct the stoichiometric/directions matrix $\Gamma$ by associating to each reaction   the column   given by $ y_{\beta}-y_{\alpha}$.

    \eeR

    The RN representation \eqr{RNSIR} may be visualized via a graph having the complexes as vertices and reactions as edges, which is known  as the \FHJ\ (HJF) graph.
The graph may be embedded in Euclidean space by using  ``species coordinates'': ``S'' is represented  by {(1,0,0)}, ``I'' is represented  by {(0,1,0)},
    and ``R'' by {(0,0,1)}. Thus, we are associating now to each reaction both an input
    {(0,1,0)}
    and an output {(0,0,1)}, rather than just the direction {(0,$-$1,1)} (as is the case in \eqr{crnsys}). Formally, let
us  summarize this in the following definition, lifted essentially from \cite{Hong}.
    \beD \label{def:structural condi}
    \begin{enumerate}
    \item The \FHJ\ graph is the directed graph whose edges are the reactions {$y_{\alpha}\to  y_{\beta},$} and the vertices are given by the complexes $y_{\alpha}, y_{\beta}$.
    \item A \emph{linkage class} is a connected component in the \FHJ\ graph, when it is regarded as an undirected graph.
    \item A CRN is said to be \emph{weakly reversible} if each of its linkage classes is strongly connected, i.e., if there is a sequence of reactions from a complex $y$ to another complex $y'$, then there must exist a sequence of reactions from $y'$ to $y$.
    \item The \emph{stoichiometric subspace} for a CRN is defined as
    $$S:={\rm{span}}\{y_r'-y_r:y_r \to y_r' \in \mathcal{R}\}.$$

    The vector $y_r'-y_r$ associated with a reaction $y_r\to y_r'$ is called a \emph{stoichiometric vector}. The matrix whose column vectors are the stoichiometric vectors is called a \emph{stoichiometric matrix}.
    \item  {For $a \in \mathbb{R}^d_{>0}$, the \emph{stoichiometric compatibility class} containing $a$ is $a+S$.} In other words, the stoichiometric compatibility class containing $a$ is the maximal set that can be reached by the deterministic system which starts from $a$.
    \item The \emph{deficiency} of a CRN is defined as $\delta:=|\mC|-l-s$, where $|\mC|$, $l$, and $s$ are the number of complexes, the number of linkage classes, and the dimension of the stoichiometric subspace, respectively.
    \item The \emph{order} of a reaction $y_r \to y_r'$ is {$y_{r,1} + \cdots + y_{r,d}$}. Reactions of orders one and two are called \emph{monomolecular} and \emph{bimolecular} reactions, respectively.
    \item The ODE associated to a CRN is
 \begin{align}\label{eq:det system}
     \frac{d{x}(t)}{dt}=\sum_r\mathcal K_r({x}(t))(y_r'-y_r),
 \end{align}
 where  $\mathcal K_r: \mathbb R^d_{\ge 0} \to \mathbb R_{\ge 0}$ is a \emph{rate function} which indicates the rate of the reaction $y_r \to y_r'$. In the mass-action case, this becomes
\be{RN} \frac{d{x}(t)}{dt}=\Gamma diag({\bff k})x^{Y_{\alpha}},\ee
where $Y_{\alpha}$ is the matrix whose columns are the source vectors $y_{r}$, and ${\bff k}$ is the vector of constants of each reaction.

\end{enumerate}
\eeD

\beXa
For  example, consider a  SIRS model with demography, with eight reactions, defined by the reactions \{0 {$\to$} ``S'',``S''+``I''  {$\to$}  2 ``I'',``I''  {$\to$}  ``R'',``R''  {$\to$}  ``S'',``S''  {$\to$}  ``R'',``S''  {$\to$}  0,``I''  {$\to$}  0,``R'' {$\to$} 0\}, or
\be{SIRS}
\frac{d{x}(t)}{dt} =
\left(
\begin{array}{cccccccc}
 1 & -1 & 0 & 1 & -1 & -1 & 0 & 0 \\
 0 & 1 & -1 & 0 & 0 & 0 & -1 & 0 \\
 0 & 0 & 1 & -1 & 1 & 0 & 0 & -1 \\
\end{array}
\right) \left(
\begin{array}{c}
 \lambda  \\
 \beta  i s \\
 i \gamma _i \\
 r \gamma _r \\
 s \gamma _s \\
 \mu  s \\
 i \mu _i \\
 \mu  r \\
\end{array}
\right).
\ee

 It has two linkage classes (see Figure \ref{f:SIRS}), one involving the  complexes $(S,I,R,0)$, and the other involving $(S+I, 2I)$. The SM has rank 3; hence, the deficiency is 6 $-$ 3 $-$ 2 = 1.

\begin{figure}[H]
\includegraphics[width=13cm]{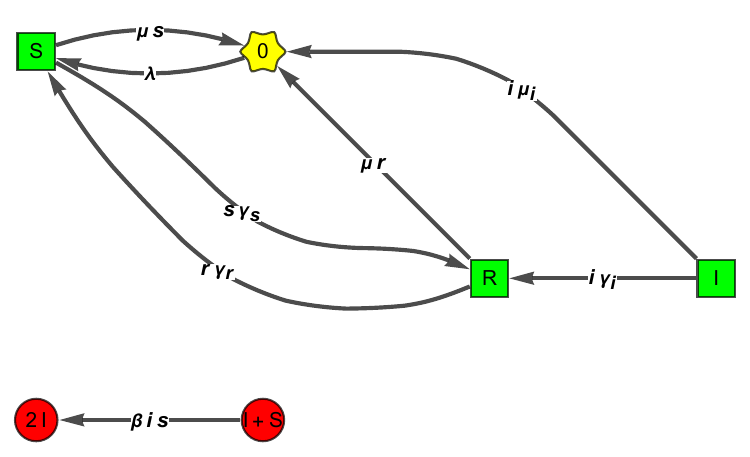}
\caption{The FHJ graph of the SIRS  with demography \eqr{SIRS} renders clear that the CRN is  not weakly reversible. The orders of the eight reactions are (0,2,1,1,1,1,1,1).}
\label{f:SIRS}
\end{figure}

SIRS   \eqr{SIRS} has two fixed points:
\BEN \im $(s= \frac{\lambda  \left(\mu +\gamma _r\right)}{\mu  \left(\mu +\gamma _r+\gamma _s\right)},i= 0,r= \frac{\lambda  \gamma _s}{\mu  \left(\mu +\gamma _r+\gamma _s\right)}$, (DFE)
\im $(s= \frac{\gamma _i+\mu _i}{\beta },i= \frac{\beta  \lambda  \left(\mu +\gamma _r\right)-\mu  \left(\gamma _i+\mu _i\right) \left(\mu +\gamma _r+\gamma _s\right)}{\beta  \mu  \gamma _i+\beta  \mu _i \left(\mu +\gamma _r\right)},r= \frac{\beta  \lambda  \gamma _i-\left(\gamma _i+\mu _i\right) \left(\mu  \gamma _i-\mu _i \gamma _s\right)}{\beta  \mu  \gamma _i+\beta  \mu _i \left(\mu +\gamma _r\right)})$ (E).
\EEN
\eeXa

\beXa [An associated reversible monomolecular model with deficiency zero] \lab{e:linSIRS}
For comparison, consider now a ``monomolecular version'' obtained by ``simplifying/translating'' the second reaction ``S''+``I'' {$\to$} 2 ``I'' to ``S'' {$\to$} ``I'' (we will say more about translation below). (see {Section} \ref{s:Ton}, and Remark \ref{r:tow} 
).

\begin{figure}[H]
\includegraphics[width=10cm]{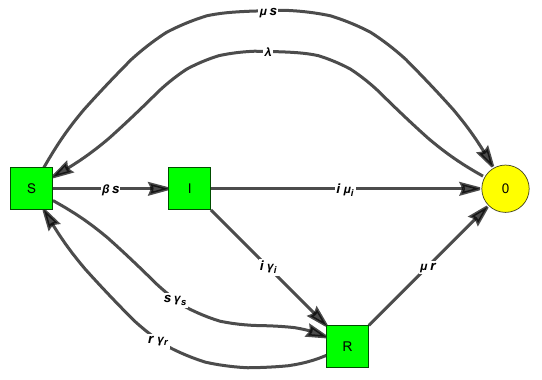}
\caption{{The}
  ``monomolecular SIRS'' has one linkage class, is   weakly reversible (WR), and has zero deficiency (ZD) 4 $-$ 3 $-$ 1 = 0.}
\label{SIRSlin}
\end{figure}

This process  has a unique globally attractive fixed point, which is positive for all values of the parameters. In fact, this is conjectured to hold always for the class of WR-ZD CRNs by the so-called ``Global Attractor
Conjecture'' \cite{Horn,CF05,CF06}. However,  this simplified   model does not make  sense epidemiologically.
\eeXa

\beR We have hinted in the last two examples one possible interest  of the CRN-ME collaboration. ME models display often multistability and oscillations (which are always proved by example dependent methods). On the other hand, they  have typically   ``WR-ZD cousins'', obtained by translation,  which are super well behaved, and which may be obtained by off-shelf software (involving sometimes Linear Programming) \cite{Hong}.

It is conceivable that the well-behaved cousins can help in reaching interesting conclusions for their more complex relatives, but this remains to be seen.  For now, we only know that product form  stationary distributions for certain CTMC stochastic CRNs (see the Section \ref{s:CTMC}) may be obtained this way \cite{HongSon}.
\eeR



\subsection{Three Crucial Differences  {Between} Mathematical Epidemiology and Chemical Reaction Network  Models}
  The most fundamental aspect  of mathematical epidemiology is the existence of at least two possible special fixed points, a boundary one, and an interior one.
   The first, the  DFE,  corresponds to the elimination of all compartments involving sickness. The second, called the endemic point (E), is  an interior fixed point which takes the ``stability relay'' from the DFE when this becomes unstable.

   This dichotomy between the stability of the boundary or that of the interior is less frequent in CRN, though it is encountered also
   in ACR (absolute concentration robustness) models \cite{AAHJ}.

   \Fr   CRN  and ME models lie often on opposite sides of the ``weak-reversibility boundary'', rendering the intersection of
 these two sets small.

 As a  last  difference, let us mention that  CRN investigates mostly models with conservations. While this is also true in a large part of   the  ME literature, which assumes a constant population, it is not true for ``varying population ME models'', which do not neglect the amount of deaths, and are therefore more challenging.


\subsection{First Fundamentals of CRNT: Complex Balancedness, Deficiency Theorems, and Toricity Conditions---CoxTcell.nb, CoxWegscheider.nb, CoxFeinberg.nb}

For a delightful introduction to CRN, the reader might want to consult, in parallel with reading this section, the slides (\cite{CoxSl}, p. 209)---see also their guide \cite{Cox}, and the experiment with the notebooks indicated in the title.

We will attempt our own brief tour of   CRN theory by recalling the celebrated  zero deficiency  (ZD) theorem of Horn and Jackson \cite{horn1972}, and of
  Feinberg \cite{feinberg1987,feinberg2019foundations} :
  \beP [ZD theorem]
 \textit{Regardless of the choice of parameters $k_r$}, a
reaction network $\mN$ with deterministic mass-action kinetics (MAK) that is both
a) weakly reversible (WR) and b) has a deficiency of zero will \saty\ the following:
\BEN
\im $\mN$ has  a ``robust ratio'' in each pair of complexes $y$ and $y'$ belonging to a common linkage class, i.e., the ratio $\mathbf{x}^{y}/\mathbf{x}^{y'}$ takes the same value at every positive steady state $\mathbf{x} \in \mathbb{R}_{> 0 }^n$.
\Mp\ \be{tre}
\frac{\mathbf{x}^{y}}{\mathbf{x}^{y'}} = \frac{{K}({y})}{{K}({y}')},\ee
where the ``tree constants'' ${K}({y})$ can be computed via a graphical procedure called the matrix tree theorem---see, for example, \cite{CoxSl}.

\im $\mN$ admits a
``complex-balanced''  equilibrium, and all its equilibria are complex-balanced.
\EEN
\eeP

 The complex balancedness from the second part of the theorem is a generalization of the well-known concept of reversible equilibrium from statistical physics, and parallel the concept of the partially balanced stationary state of probability, which we will take for granted here. What must be recalled for our paper,  and for the development of computer packages, is  its Definition \ref{d:CB} below, which requires a  matrix representation of mass-action CRNs due to \cite{Horn}
\be{Lap}  x'=\Gamma {\bf R}(x)=(\Yb I_E ) (I_k x^{\Yb})=\Yb \mL  x^{\Yb}.\ee

 Here, $\Yb$ is a $n_S\times n_C$ matrix whose columns specify the composition of the ``complexes'' (i.e., the vectors representing the
left- and right-hand sides of each equation); $I_E $ is the $n_C\times n_R$ incidence matrix of the \FHJ\  graph on the complexes, whose edges correspond to the reactions; $x^{\Yb}$ is the  $n_C\times 1$ vector of exponents  of the inputs of each reaction, completed by $1$ for the complexes which are not inputs;
    and  $I_k$ is an $n_R\times n_C$ ``outgoing coincidence matrix'' containing all the reaction constants, whose $(r, c)$-th element equals $k_r$ if the $c$-th complex is the input complex for the j-th reaction, and zero
otherwise.

It follows that $\mL=I_E  I_k$ is a $n_C\times n_C$ ``Laplacian matrix'' with non-positive diagonal elements and non-negative off-diagonal
elements, whose column sums  are  zero. This matrix intervenes in the study of ``complex-balanced equilibria'' and Wegscheider conditions---see, for example \cite{Gater,van2015complex} for further details.

\beD \lab{d:CB}
 A mass-action CRN   defined by  \eqr{Lap} is called complex-balanced (CB) at an equilibrium point $x_* \in \mathbb{R}_+^n$ if $$I_E  {\bf R}(x_*)=\mL x_*^{\Yb}=0.$$

\eeD

\beR 

Note the representation \eqr{Lap} suggests a conceptual decomposition of CRN ODEs into (1) the \sm\ (SM); (2) a weighted graph structure specified by a Laplacian matrix; and (3) monomial rates.

This gives rise to the further generalization  of \GMAK\ (GMAK) \cite{Horn,MulReg}, in which the last matrix in \eqr{Lap}, to be called the kynetic matrix (KM), may be different from the first matrix (SM).

\eeR

\beR If a polynomial dynamical system admits a deficiency zero realization that is {\em not} weakly reversible, then its dynamics is also greatly restricted: it can have no positive steady states, no oscillations, and no chaotic dynamics~\cite{Horn,feinberg2019foundations},
 (\cite{feinberg1987}, Remark 6.1.B) for any choice of parameters $\kappa_r$.

This shows that such models cannot appear in ME (because  we are only interested in models that admit a positive steady state when $R_0>1$).
\eeR
\Oth models with both ZD and WR also probably cannot appear in ME, due to the ``Global Attractor
Conjecture'',   which states that the complex-balanced
equilibria of reaction networks are globally asymptotically stable
relative to the interior of their positive compatibility classes \cite{Horn}. Later, this  conjecture, whose name was given by Craciun {et al.} 
 \cite{CraStu}, became the Graal of CRN;  in the case of single-linkage CRNs, it was proved by Anderson
 \cite{anderson2011proof}.

Thus, deficiency zero is probably impossible for ME models; however, a few simple models   like SIRS   \eqr{SIRS} do have deficiency one but do not \saty\ the conditions of the deficiency one theorem.

The above  suggests that immediate applications of  CRN methods to ME models are not easy  to find---see, however, the next subsection for some recent ``intersection results''.


We provide now some  examples from \cite{Cox}, implemented in the notebooks mentioned in the title.  The first illustrates the power of the ZD theorem.
\beXa The system
$$ \fr d{dt} \bep a\\b\\c\\d \eep=\bep-a b k_1+c k_2+d k_4\\-a b k_1+c k_2+d k_4\\a b k_1-c \left(k_2+k_3\right)\\c k_3-d k_4\eep$$
is quite a challenge analytically since it may be checked to have a 2-dim set of fixed points, parameterized by a,c, with a complicated singularity
in the origin.  However, the ZD theorem allows us to conclude directly that  {{for every set}} of
positive reaction rates, in any stoichiometric compatibility class,
there is a unique positive steady-state solution which  is complex-balanced and locally attracting, and oscillation is impossible.

\eeXa

While ZD-WR CRNs are complex-balanced (also called toric) for every choice of positive constants, certain networks, as known since Wegscheider~\cite{Weg}, are only CB/toric under certain conditions.  Toricity conditions are beyond the scope of this review---see \cite{CoxSl,Cox}.  We will offer though some  notebooks where they are determined by Mathematica's elimination tools.
\beXa The system (\cite{CoxSl}, p. 241)
$$ \fr d{dt} \bep a\\b \eep=\bep-\left(a^2 \left(\kappa _{31}+2 \kappa _{32}\right)\right)+a b \left(\kappa _{13}-\kappa _{12}\right)+b^2 \left(2 \kappa _{23}+k_{21}\right)\\a^2 \left(\kappa _{31}+2 \kappa _{32}\right)+a b \left(\kappa _{12}-\kappa _{13}\right)-b^2 \left(2 \kappa _{23}+k_{21}\right)\eep$$
is used as an illustration of the matrix tree theorem, which yields the Wegscheider-type  toricity condition $K_1 K_3 =K_2^2$, where $K_1=\kappa _{21} \kappa _{31}+\kappa _{23} \kappa _{31}+\kappa _{21} \kappa _{32}, K_3=\kappa _{13} \kappa _{21}+\kappa _{12} \kappa _{23}+\kappa _{13} \kappa _{23},K_2=\kappa _{12}\kappa _{32}+\kappa _{13}\kappa _{32}+\kappa _{31}\kappa _{12}$.

\eeXa

\beXa The Edelstein--Feinberg--Gatermann example (\cite{CoxSl}, p. 252)
$$ \fr d{dt} \bep a\\b \\c\eep=\left(
\begin{array}{c}
 a \left(-a k_{21}-b \kappa _{34}+\kappa _{12}\right)+c k_{43} \\
 c \left(\kappa _{45}+k_{43}\right)-b \left(a \kappa _{34}+\kappa _{54}\right) \\
 a b \kappa _{34}+b \kappa _{54}-c \left(\kappa _{45}+k_{43}\right) \\
\end{array}
\right)$$
has toricity condition $\kappa _{12}=\frac{\kappa _{54} k_{21} k_{43}}{\kappa _{34} \kappa _{45}}$.

\eeXa

Here is a stochastic counterpart of the deficiency zero theorem,  provided by
Anderson et al. \cite{Anderson2010}.
\beP
To a deterministic CRN $(\mS,\mC,\mR)$,  associate a CTMC

$X(t)=(X_{[1]}(t),X_{[2]}(t),\dots,X_{[d]}(t))\in Z^d_{\ge 0}$ whose transition rates are given by
\begin{align*}
    P(X(t+\Delta t)=\mathbf{n}+\zeta \ | \ X(t)=\mathbf{n})=\sum_{k:y_r'-y_r=\zeta} \lambda_r(\mathbf{n})\Delta t + o(\Delta t) \quad \text{for each $\zeta \in \mathbb Z^d$,}
\end{align*}
where the $i$-th coordinate  represents the  number of species $S_i$, and $\lambda_r:\mathbb Z^d_{\ge 0}= \mathbb R_{\ge 0}$, the \textit{mass-action propensity function} associated with the reaction $y_r \to y_r'$, is given by
\begin{align}\label{eq: stoch mass}
   \lambda_r(\mathbf{n})=\kappa_r\dfrac{\mathbf{n} !}{(\mathbf{n}-y_r)!}\mathbf{1}_{\mathbf{n} \ge  y_r},
\end{align}
and where we note that the powers in the continuous mass-action kinetics
   $  \mathcal K_r(x)=\kappa_rx^{y_r}
 $ have been replaced by decreasing factorials.

Now let $(c_{[1]}, \ldots, c_{[d]})$ denote  the complex-balanced equilibrium of a deterministic  CRN which has WR and ZD. Then, the stochastic  CRN under  stochastic mass-action kinetics \eqref{eq: stoch mass} admits a Poisson product-form stationary distribution that is given by
\[
 \pi(\mathbf{n})= M\displaystyle \prod_{i=1}^d \frac{c_{[i]}^{n_{[i]}}}{n_{[i]}!},
\]
where  $M>0$ is a normalizing constant.
\label{anderson:2010:theorem}
\eeP

We conclude this section by noting that the remarkable simplifications for WR and ZD networks gave rise to the natural question of whether one might find CRN realizations whose associated dynamical system is equivalent to a given one, and also has one of these two properties. One such method,
called ``Network Translation'',  initiated by Johnston~\cite{John,tonello2018network}, gave rise to several applications and software implementations, for example, TOWARDZ \cite{Hong}.  For example, (\cite{tonello2018network}, {Thm 2,3}) showed that formula \eqr{tre} 
holds also whenever a MAK may be translated into a ZD-WR GMAK.
A further example illustrating the combined power of the ZD theorem and of network translation will be given in the last Sections \ref{s:Ton} and  \ref{s:tow}.

Further applications for CTMC CRNs are provided in \cite{HongSon,hoessly2023complex}.

\subsection{Some Recent Interactions  Between CRN and ME Methods}\lab{s:syn}
\BEN \im The idea behind the representation \eqr{SIRS0}, see also Remark \ref{r:str}, was exploited for ME models in the recent paper \cite{AAV}. An interesting result there is (\cite{AAV}, Thm 3.1), which  may be informally stated as follows:
      Consider any epidemiological system in which there exists an ``S{$\to$}I infection  reaction'' with admissible rate $R_1(s,i)$, and an ``I{$\to$}\dots  treatment  reaction'' with  admissible rate $R_2(i)$. Then, the ``symbolic Jacobian'', in which $R_1,R_2$ are not specified (but the sign of their derivatives
      is specified, via admissibility), may always have purely imaginary eigenvalues if $R_1,R_2$ are ``rich'' in the sense of \cite{Vas,VasStad} (for example, \MM).

      Briefly,
 all  epidemic models admit ``symbolic  bifurcations'', provided their rate functions  have enough parameters;  it is only the restriction
to mass action that may prevent the occurrence of Hopf bifurcations.

Note that this result had been observed  empirically in many particular three-compartment SIR-type models, but the fact that the  number of  compartments and the exact architecture of the model are irrelevant was not properly understood.

 \im A second CRN result  exploited in \cite{AAV} was the {\bf {Inheritance of oscillation in chemical reaction networks}} of  \cite{banaji2018inheritance}, which gives conditions for Hopf  bifurcations to be inherited by models, given that they exist in a submodel where some parameters are $0$.
  This allowed establishing the existence of  Hopf  bifurcations for a certain mass-action ME model, using their existence in a simpler case already studied by Hethcote and Van den Driesche \cite{HethSVDD}).
  \im In a recent preprint \cite{AAHJ}, we provided an interaction in the opposite sense, by showing  that the stability of boundary points for the class of CRNs with \ACR\ (ACR) may be analyzed via the NGM method of ME.
\EEN

 \beO \lab{o:MEm}
 Mathematical epidemiology ODE models  could be defined (hopefully with benefits) as  particular CRN models  formed with  only three types of reactions:
\BEN \im Transfers (monomolecular reactions);
\im  Bimolecular auto-catalytic reactions of the type $S+I\to 2 I$  as encountered in SIR, etc;
\im  Bimolecular auto-catalytic reactions of the type $ S+I\xrightarrow{[\beta_e s i]}I+E \Eq s'=-\beta_e s i \dots , e'=\beta_e s i\dots $ as encountered in SEIR.
\EEN

In addition, they should have at least one boundary fixed point and one interior fixed point.
\eeO
\beR By the suggested definition  above, the EnvZ-OmpR model in \eqr{envzompr} below cannot be accepted as an ME model, even though its qualitative behavior  is similar to that of ME models. \eeR

\beO Can an ME model defined as in Open Problem \ref{o:MEm} be weakly reversible?\eeO

\beR As an aside, for  mathematical virology ODE models, it seems that most of the bimolecular reactions encountered in CRN (with sum of the coefficients less than two in  the LHS  of each reaction) might be of practical relevance.
\eeR

\subsection{Can CRN Software Solve ME Problems?}
 CRN methods were implemented
 in  powerful  software, like, for example, the following: \BEN \im
 The   collaborative package {CoNtRoL} (accessed on 6 August 2015) 
 \cite{donnell2014control,Control} \url{https://reaction-networks.net/wiki/CoNtRol}{;} 
 \im The  Mathematica packages ReactionKinetics \cite{Toth}, {MathCRN}  (accessed on 29 April 2016)
 \url{https://gitlab.com/csb.ethz/MathCrn} and
 reaction-networks.{m} \url{https://vcp.med.harvard.edu/software.html\#crnt};
 \im Feinberg's  Chemical Reaction networks {toolbox} \url{https://cbe.osu.edu/chemical-reaction-network-theory#toolbox};

 \im The Julia package Catalyst.jl \cite{catalyst};
 \im The Macaulay package ReactionNetworks.{m2} \url{https://macaulay2.com/doc/Macaulay2/share/doc/Macaulay2/ReactionNetworks/html/index.html};
  \im The  Matlab package {LEARN} \cite{AAnew,AASComp,
Ali23gr} \url{https://github.com/malirdwi/LEARN}{;} 
\im The Matlab packages TowardZ and CASTANET for network translation, a topic initiated by {Johnston}  \cite{John,tonello2018network,JMP,hernandez2022framework}), implemented by {Hong et al.}  \url{ https://github.com/Mathbiomed/TOWARDZ} {and}
 \url{https://github.com/Mathbiomed/CASTANET};
\im  The Python packages \cite{Tonello,crnpy} {and}
 \url{https://github.com/materialsproject/reaction-network}.
\EEN

See {also}
 \url{https://reaction-networks.net/wiki/Mathematics_of_Reaction_Networks} for further resources.

It is interesting to reflect on  the reason for the disparity between the strong computerization of CRN and its absence in ME for similarly looking problems. One reason is that ME strives towards complete analysis of small-dimension models (but see the 31 reactions model of \cite{Bulh} for an exception), {while CRN researchers study larger models but  settle  for partial answers (concerning, say, multistationarity and the existence of ACR) and avoid often difficult issues like the presence of Hopf and Bogdanov--Takens bifurcations, or chaos.}

\subsection{Contents}

Section \ref{s:ME}  {revisits some key concepts  in mathematical epidemiology, among them being the stability of boundary points}, which  are illustrated by  two examples, the  SAIR and SLAIR epidemic models.

Section \ref{s:Ton} touches briefly on some CRN state-of-the-art topics
which might turn out to be useful for ME researchers in the future.

Section \ref{s:SIRCRN}  concludes with a   discussion of future perspectives. 

 {Section} 
 \ref{s:CTMC} reviews    some relevant facts on continuous-time Markov chain  (CTMC) models, which lurk in the backstage since ODE models may be viewed as limits of CTMCs.

\section{A Revisit of Some Key Concepts in Mathematical Epidemiology: The Disease-Free
Equilibrium, the Next Generation Matrix,  the Basic Reproduction Number \boldmath$R_0$, and the Diekmann Kernel \label{s:ME}}

We will start our tour of ME by an example of ME model, a generalization of the classic SEIR model.

\subsection{An Example: The SAIR/S$I^2$R/SEIR Epidemic Model\label{s:SEIR}: SAIRS.nb}

We introduce here a {nine}-parameter SAIR/S$I^2$R/SEIR-FA epidemic model \cite{Van08,RobSti,Ansumali,Ott,AAH} as a concrete example for the generalization
presented in the next section, which is parameterized by two matrices.

 We now warn the reader that the equations we choose to study below, for the fractions  $\s, \mathsf a, \mathsf i$, and $\mathsf r$  of susceptible, asymptomatic, infected, and recovered, are an only an approximation (albeit a very popular one).  Indeed, the model defined  in Figure \ref{f:SAIRS} has
 varying population models due to $\de>0$, and the correct model for the fractions must include further quadratic terms multiplied by $\de$ \cite{LiGraef,AABH,AAH}:
 \begin{equation}\label{SEIRsc} \bc
\s'(t)= \Lambda   -\s(t)\pr{\beta_i  \mathsf i(t)+\beta_a \mathsf a(t)+\g_s+\Lambda} + \g_r \mathsf r(t)\\
\bep \mathsf a'(t)\\\mathsf i'(t)\eep =
\pp{\s(t) \bep  \beta_a     &  {\beta_i} \\
 0&  0 \eep- \bep    \gamma_a +\Lambda  &  0 \\
 -\ei& {\Lambda + \g_i + \de}\eep}
  \bep \mathsf a(t)\\\mathsf i(t)\eep
\\
\mathsf r'(t)=  \g_s \s(t)+ a_r   \mathsf a(t)+ {\g_i}   \mathsf i(t)- (\g_r+\Lambda) \mathsf r(t) \ec,
\end{equation}

\begin{figure}[H]
 \includegraphics[width=13.5cm]{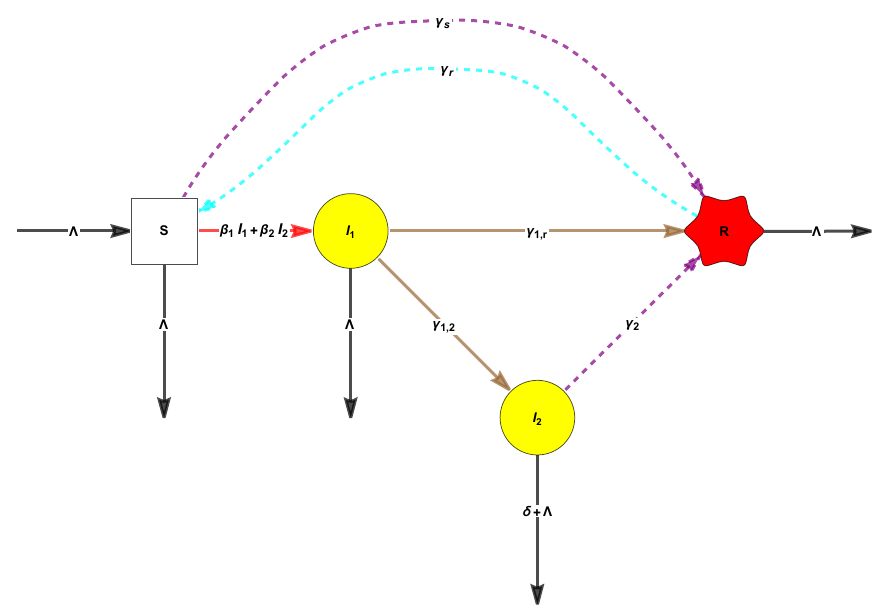}
 \caption{Chart flow/species graph of the SI$^2$R model, with two infected classes and extra deaths at rate $\delta$.  The red edge corresponds
to the entrance of susceptibles into the disease classes, the brown edges are the rate of the transition matrix V, and the cyan
dashed lines correspond to the rate of loss of immunity.  The remaining black lines correspond to the
inputs and outputs of the birth and natural death rates, respectively, which are equal in this case.}
 \label{f:SAIRS}
\end{figure}

\beR~
\BEN
 \im We write the ``infection'' middle equations in the form $\vi'= M \vi$, where $\vi=\bep \mathsf a(t)\\\mathsf i(t)\eep$ to emphasize their factorization.
Also, for the factor appearing in these equations, we emphasize a form
 \begin{equation}\label{V}F - V.\end{equation}

 Such decompositions, not necessarily existing nor unique \cite{Diek,Van,Van08}, are  used in the computation of the  next generation matrix (NGM)
 \be{K}K=F.V^{-1}\ee
 and of its spectral radius, the \brn\ $R_0$.
 \im  The  SAIR model is obtained when $a_r{=}(\g_{1,r})=0=\de$, and the classic SEIR model is obtained when \frt\ $\beta_{a}=0$.
\EEN \eeR

The reactions of the corresponding mass-action CRN are

$0 \xrightarrow{\Lambda} S,
S+A \xrightarrow{a s \beta _a} 2 A,
 {S+I \xrightarrow{i s \beta _i} A + I},
 {A \xrightarrow{ a \left(\gamma _a-a_r\right)} I},
 {A\xrightarrow{a \left(\gamma _a-a_i\right)} R},
 {S \xrightarrow{ s \gamma _s} R},\\
 {I \xrightarrow{ i \gamma _i} R},
 {R \xrightarrow{ r \gamma _r} S},
 {S \xrightarrow{ \Lambda  s } 0},
 {A \xrightarrow{  \Lambda a } 0},
 {I \xrightarrow{ (\Lambda+\delta) i } 0},
 R \xrightarrow{  \Lambda r} 0,
 $

 see the Figure \ref{SAIRgr} below.

\begin{figure}[H]
\includegraphics[width=13cm]{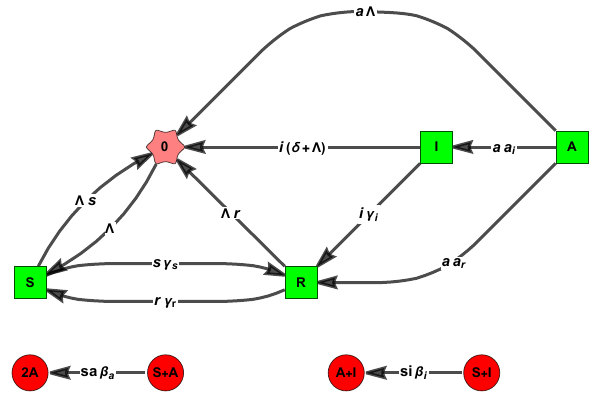}
\caption{{The} 
 {Feinberg--Horn--Jackson graph} of the ``SAIR network'' with $n_V=9$ vertices $(S,A,I,R,S+A,S+I, 2A, A+I,0)$ (where the 0 node represents the exterior), 12 edges, and 3 linkage classes.
The deficiency  is $n_V-rank(\Gamma)-n_C=9-4-3=2$, and weak reversibility does not occur, so neither of the conditions for having complex-balanced equilibria holds.}
\label{SAIRgr}
\end{figure}



\subsection{{The} Next Generation Matrix Approach}

In this section, we discuss the computation of the stability domain of the DFE, via the fascinating  \NGM\ (NGM) method \cite{Diek,Van,Van08}, called this way since it replaces the investigation of the Jacobian with that of a matrix whose origins lie in probability (the theory of branching processes).  The NGM approach seems to  be very little known outside ME,  but we have shown in \cite{AAHJ} that it may be applied to
CRN models with  ACR (absolute concentration  robustness).

The NGM method is based on the projection of the Jacobian on the subset of infection variables, and is justified by conditions which may be found in
\cite{Diek,Van,Van08}. These seem to be satisfied in all ME models with polynomial rates, so instead of stating them, we will offer here a less known ``NGM heuristic'', consisting of the following steps:
\BEN \im Write the infection equations (the middle ones in the Section \ref{s:SEIR}) in the form
\be{M}\vi'= M \vi\ee
\im Split $M$ into the part of ``new infections'' $F$
containing  all the non-constant terms with positive sign \cite{AABJ}, and  the rest, denoted by $-V$, arriving thus to
\be{FV} \vi'= (F-V) \vi\ee
The purpose of this is to replace the study of the spectrum of Jacobian $M$
with that of the ``next generation matrix'' defined by
\be{K} K=F.V^{-1}.\ee
\im Determine the spectral radius of $K$, the so called \brn\ $R_0$.\EEN

\beR The first step is suggested by the well-known fact that infection variables are fast compared  to the others, so the non-infection variables may be taken as fixed and their equations ignored, asymptotically.

The second  step is justified by an empirical observation made in \cite{AABJ}. We recall that  the rigorous NGM is obtained via a decomposition which was not specified uniquely in the  papers cited above. This leads to the non-uniqueness of $R_0$, and leaves its final choice to the latitude of the ``expert epidemiologist'', a situation which is maybe not ideal.
We propose in \cite{AABJ} to complement the classical NGM method by the unique specification of   $F$ described above, and show that this produces reasonable answers in all the examples we investigated.

Let us mention also another empirical observation, that the \ch\ of the matrix $K$ factors usually more than the \ch\ of the Jacobian of the infection equations \wrt\
the infection variables, which explains the popularity of the NGM method.

\eeR

For an example of applying the NGM method, see Remark \ref{r:hNGM}.

In the next section, we  present a special but quite general class,  the SIR-PH-FA  models introduced by us in \cite{AABBGH}, which  generalize the SAIR example from Section \ref{s:SEIR}, and which we deem important for several reasons, to be explained at the end of  Section \ref{s:ker}.

We conclude by explaining the  acronym SIR-PH-FA. The PH (phase type), is due to the probabilistic interpretation of SIR models where the passage through the infection classes, before recovery, is similar to the lifetime of an absorbing CTMC with generating matrix 
$(-V)$, which has a phase-type distribution. The FA (first approximation) part is due to the fact that these models will have, in general, a varying population, and that it is convenient to study the corresponding proportions  from the total population. Since these are quite hard to study in general---see, for example, \cite{LiGraef,AABBGH}---it is common practice to ignore the number of deaths (quite acceptable for short periods), and the first example of this is the SIR model of \cite{Ker}. This practice is so common that, in fact, all ME researchers refer to FA ME models as ME models, simply. We cannot do that, however, because  we have introduced in \cite{AABBGH}
a more refined IA (intermediate approximation of varying population models), and we want to distinguish between the extensively studied FA models, and our own (very little studied) IA models.

\subsection{SIR-PH-FA Models} \lab{s:renker}

We revisit here a class of     epidemic models \cite{AABBGH,AABH,AABGH}, which is a  particular case of the yet more general \AB\  introduced in \cite{AABBGH}. Here, PH stands for ``phase-type'' (distribution), and FA (first approximation) stands for the fact that certain terms from the correct SIR-PH model have been neglected, which will be further explained below. 

Before proceeding, we need to comment on the ``intrusion'' of PH, yet another Markovian concept (the first two where the NGM and the \brn) into the theory of ODE models, and also to clarify that this third intruder is not sufficiently known since, unlike the first two, it was only introduced recently in a paper that was never published, by Riano \cite{Riano}.

\beR \lab{r:PH} To clarify the above statement further, the fact that ODE models and their ``CTMC versions'' are very close to each other  is very well known (and reviewed in the Section \ref{s:CTMC}). 
What was very little mentioned in the literature prior to \cite{Riano} was the fact that the ODE model projected on the infection equations
also has as a close version, an absorbing CTMC. More precisely, the matrix $-V$ in \eqr{FV} may be viewed as the generator of a Markovian evolution
among the infection states, prior to their absorption in the recovered state (and the susceptible state, if second infections are possible).
As one echo of this possibly forgotten fact, note that the elements of $V^{-1}$ in \eqr{K} are sometimes called  "expected dwell times'', adopting
the terminology of absorbing CTMCs---see {Subsection} \ref{s:PH}.
\eeR

SIR-PH-FA    epidemic models are parameterized essentially by two matrices, which we introduce now via the SAIR example \eqr{SEIRsc}:
\bea&& \bc
\s'(t)= \Lambda   -\s(t)\pr{\beta_i  \mathsf i(t)+\beta_a \mathsf a(t)+\g_s+\Lambda} + \g_r \mathsf r(t)\\
\bep \mathsf a'(t)\\\mathsf i'(t)\eep =
\pp{\s(t) \bep  \beta_a     &  {\beta_i} \\
 0&  0 \eep- \bep    \gamma_a +\Lambda  &  0 \\
 -\ei& {\Lambda + \g_i + \de}\eep}
  \bep \mathsf a(t)\\\mathsf i(t)\eep
\\
\mathsf r'(t)=  \g_s \s(t)+ a_r   \mathsf a(t)+ {\g_i}   \mathsf i(t)- (\g_r+\Lambda) \mathsf r(t) \ec
\eea
and let us rewrite the middle equations in the form:
\bea &&\bep \mathsf a'(t)\\\mathsf i'(t)\eep =
\pp{\s(t) B +A - \bep    \Lambda  &  0 \\
 0& {\Lambda  + \de}\eep}
  \bep \mathsf a(t)\\\mathsf i(t)\eep \Eq M=\s(t) B +(A - \bep    \Lambda  &  0 \\
 0& {\Lambda  + \de}\eep)
 \no\\&&=F-V, F=\s(t) B, V=-A + \bep    \Lambda  &  0 \\
 0& {\Lambda  + \de}\eep,\eea
where $B=\bep  \beta_a     &  {\beta_i} \\
 0&  0 \eep, A=\bep   - \gamma_a   &  0 \\
 \ei& {- \g_i }\eep$.

 {Here, the constant matrix $B$ gathers all the infection rates,  while $A$ is the generator of the phase-type  semigroup describing the CTMC, which may be associated with transitions between the infected compartments, prior to absorption in the susceptible or recovered classes----recall Remark \ref{r:PH}.
 Finally, the diagonal matrix containing demography terms has the probabilistic  interpretation of extra killing (moving outside the four compartments).

 \beR The  matrix $B$ has rank one due to the existence  of only one susceptible class, and  this simplifying feature will be kept in the generalization presented next. \eeR

{We recall now the powerful generalization of SAIR in \cite{AABBGH}, suggested by the works of~\cite{Arino,Riano}.}

\beD \lab{d:SIRPH} A SIR-PH-FA   model (\cite{AABH}, Definition 1) is defined by
\begin{align}
\label{SYRPH}
\vi'(t) &=  \vi (t) \pp{\s(t)  \; B+ A
- Diag\pr{\bff{\de}+\La \bff 1}}:=
 \vi (t) \boxed{(\s(t) B-  V )}  \no\\
 \s'(t) &= \pp{\La -\pr{ \La + \g_s } \s(t)}- \s(t) \T {i}(t)+ \g_r {\mathsf r}(t), \; \; \T {i}(t) = \vi(t) \bb \nonumber\\
 \bb&=\bep\bb_1\\\vdots\\\bb_n\eep, \quad \bb_i=(B \bff 1)_i=\sum_{j} B_{i,j}, i=1,\dots ,n\nonumber\\
{\mathsf r}'(t) &=  \vi (t) \bff a+ \s(t)  \g_s   - (\g_r+\La){\mathsf r}(t), \qu \bff a=(-A) \bff 1.
 \end{align}

 Here
\BEN
\im $\s(t) \in \mathbb{R_+}$ represents the set of individuals susceptible to be infected (the beginning state).
\im ${\mathsf r}(t) \in \mathbb{R_+} $  models
recovered  individuals (the end state).
\im $\g_r$ gives  {the rate at which recovered individuals lose immunity,}
and  $\g_s$ gives the rate at which individuals are vaccinated (immunized). These two transfers connect directly the beginning and end states (or classes).
\im The row vector $\vi(t) \in \mathbb{R}^n$  represents the set of individuals in different disease states.
\im $\Lambda >0$ is the per individual death rate, and it equals also the global birth rate (this is due to the fact that this is a model for proportions).

\im  $A$,  which describes transitions between the disease classes, is a $n\times n$ {\bf Markovian subgenerator matrix}, {(i.e., with negative diagonal elements dominating the sum of the positive off-diagonal ones).  } \Mp a Markovian subgenerator
matrix satisfies that each
off-diagonal entry  $A_{i,j}\geq 0, i\neq j$,  and  that the row-sums $\sum_j A_{i,j} \le 0,
\for i$  with  at least one inequality being strict.

{Alternatively, $-A$ is a  non-singular M-matrix \cite{Arino}, i.e.,
  a real matrix $P$ with  $ P_{ij} \leq 0, \forall i \neq j,$ and having eigenvalues
whose real parts are non-negative \cite{plemmons1977m}.}

The fact that a Markovian subgenerator appears  in  our ``disease equations'' suggests that certain probabilistic concepts intervene in our deterministic models, and this is indeed the case---see~\cite{AAB,AAK} and below. 

\im $\bff \de \in \mathbb{R_+}^n$ is a column vector giving the death rates caused by the epidemic in the  disease \com s.
\im  The matrix $-V$, which combines $A$  and the birth and death rates $\La, \bff \de$ by
\be{V} V:=-A+ Diag(\delta+ \Lambda \bff 1 )\ee
 is also a Markovian subgenerator. This entails that
$V^{-1}$ contains only positive elements, which are precisely the expected  ``dwell times'' (i.e.,  times spent in each  infection class---see \eqr{edt}) before absorption into S, in the associated CTMC model \cite{AAB}.

\im $ B $ is an $n \times n$ matrix. {We will denote by $\bb$ the vector containing the sum of the entries in each row
of $B$, namely, $\bb= B \bff1$.} Its components $ \bb_i$ represent the
{\bf total force of infection} of the  {disease} class $i$, and $\s(t)   \vi(t) \bb$ represents the
total flux which must leave class $\s$.
Finally, each  {entry} $B_{i,j}$, multiplied by $\s$,  represents the
force of infection from the  {disease} class $i$ onto class $j$, and our essential assumption below will be that $B_{i,j}=\beta_i \al_j$, i.e., that all forces of infection are distributed among the infected classes conforming to the same probability vector $\va=(\al_1,\al_2,\dots ,\al_n)$.

\EEN

\eeD

\beR \label{r:PF}

 The Jacobian of the SIR-PH-FA   model is explicit. With $\g_r=0$ (for simplicity), it~is
\begin{equation}\label{JacPh} {\bep \s B-V& \vi B \\ -\s \bb& -\La -\g_s- \vi \bb \eep}\end{equation}
where $\bb$ is defined in \eqr{SYRPH}.
\eeR

\beR
Note that the factorization of  Equation \eqr{SYRPH} for the diseased \com s $\vi$  implies a representation of $\vi$ in terms of $\s$:
\begin{equation}\label{irep}
\vi(t)=\vi(0) e^{- t V + B \int_0^t \s(\tau) d \tau} =\vi(0) e^{\pp{- t I_n + B V^{-1} \int_0^t \s(\tau) d \tau} V}. \end{equation}

{In this representation intervenes an essential character of our story, the matrix $B V^{-1},$ which is proportional for SIR-PH-FA   models to the \NGM\ $\s B V^{-1}$}. {A  second representation~\eqr{vc} below will allow us to embed our models in the interesting class of distributed delay/renewal models, in the case when
$B$ has rank one.}

\eeR

\subsection{Markovian Semi-Groups Associated to SIR-PH-FA  Epidemic Models with One Susceptible Class and \NGM\  of Rank One:  Their  ``Age of Infection
Diekmann Kernel'', and an Explicit $\mR$  Formula  for Their Replacement Number \lab{s:ker}}

The purpose of this section is to revisit, for SIR-PH-FA  epidemic models (with one susceptible class and \NGM\  of rank one), a kernel defined first in \cite{AABBGH} and below in \eqr{ker}, which generalizes conceptually the replacement number $\mR$ defined in \eqr{mR} below.

\beP\label{p:ren}  Consider a  SIR-PH-FA   model \eqr{SYRPH} with one susceptible class, with $B= \bb \va$ of rank one, and with $\g_r=0$ so that ${\mathsf r}(t)$ does not affect the rest of the system. {Let }
$$\T {i}(t) =  \vi(t) \bb$$
denote  the total force of infection.
Then, we have the following:
\BEN
\im  The solutions of the ODE system
\eqr{SYRPH}    \saty\ also a ``distributed
delay SI system''  of two scalar equations
\begin{align}\label{SI}
\bc \s'(t)=\Lambda -\pr{ \La + \g_s } \s(t)-\s(t) \T i(t)\\
\T i(t)=\vi(0)  e^{-tV} \bb+
\int_0^t s(\tau) \T i(\tau)
a(t-\tau) d\tau,\ec
\end{align}

where \begin{equation}\label{ker} a(\tau)=  \va e^{-\tau V} \bb, 
 \end{equation}
with $-V=A
- \pr{Diag\pp{\bff{\de}+\Lambda \bff 1}}$  (it may be checked that this fits the formula on page 3 of~\cite{Breda} for SEIR when $\Lambda= 0, \de=0$). ({$a(t)$ is called ``age of infection/renewal  kernel; see }\cite{Hees,Brauer05,Breda,DiekHeesBrit,Champredon,Diek18,Diek22} for expositions of this concept.)

\im Define the \bzn\ via the integral {representation} {(the}  Arino et al. \cite{Arino} formula)
\begin{equation}\label{mR}
\mR=\int_0^\I a(\tau)  d \tau= \int_0^\I \va e^{- \tau V} \bb d\tau =\va\ V^{-1} \; \bb .
\end{equation}

Then, the \brn\ and  \bzn\ \saty
\be{R0sR}R_0={\sd}\mR.\ee

\EEN
\eeP

\beR The definition \eqr{mR} comes from the ``survival method'', a first-principles method whose rich history  is described in   \cite{Hees,Diek10}---see  also (\cite{Champredon}, (2.3)), (\cite{Diek18}, (5.9)).\eeR

\begin{proof}
1. The non-homogeneous   infectious equations may be transformed into an integral equation by applying the variation of constants formula.
The first step is the solution of the homogeneous part.  Denoting this by  $\Gamma(t)$,  \ith \begin{align}
\label{SYRIGI}
\vec \Gamma'(t) &=
- \vec \Gamma(t) V \Longrightarrow \vec \Gamma(t) = \vec \Gamma(0)e^{t (-V)}.
\end{align}

{When $\vec \Gamma(0)$ is a probability vector,} \eqr{SYRIGI} has the interesting probabilistic interpretation of the survival probabilities in the various components  of the semigroup generated by the Metzler/Markovian subgenerator matrix $-V$ (which inherits this property from the phase-type generator $A$). Practically, $\vec \Gamma(t)$ will give  the expected fractions of individuals who are still in each compartment at time $t$.)

The variation of constants formula implies then that $\vi(t)$ satisfies the integral {equation} 
\begin{equation}\label{vc} \vi(t) = \boxed{\vi(0) e^{-t V} +
\int_0^t \s(\tau) \vi(\tau)  B e^{-(t-\tau) V} d \tau}.\end{equation}

Now in the rank one case $B=\bb \va$, and
 \eqr{vc} becomes
\begin{equation}\label{vcb} \vi(t) = \vi(0) e^{-t V}+
\int_0^t \s(\tau) \pp{\vi(\tau)  \bb} \va e^{-(t-\tau) V} d \tau.
\end{equation}

Finally, multiplying both sides on the right by $\bb$ yields the result.

2. See (\cite{AABBGH}, Prop. 2 and 3), or, alternatively,  note that all eigenvalues of the \NGM\ except one are $0$ \cite{Arino,AABBGH}.
\end{proof}

In conclusion, SIR-PH-FA models with \NGM\  of rank one are important due to the following:
\BEN \im  They allow relating two   approaches for computing  $R_0$, and they render evident the  importance of the \bzn\ $\mR$,  the integral of  the age of infection/Diekmann kernel (\cite{AABBGH}, (7)) (see  \cite{Diek} for the origins of this concept).

 \im  They are quite tractable  analytically due to  the existence of a unique endemic point EE, for which
 \be{sEACR} s_{E}=1/\mR\ee
 (\cite{AABBGH}, Prop. 2 and 3).  This may further be elegantly expressed
 as \be{R0fr} R_0=\fr{\sd}{s_{E}}.\ee

 \EEN

 \beR There exists some confusion in the literature between the \brn\ and the \bzn (we follow here the terminology of Hethcote \cite{Heth}).

The main points are that we need two concepts, as is evident from \eqr{R0sR}, and that $\mR$ cannot be ignored since it appears in  important formulas like \eqr{sEACR}.

 Let us mention  that while \eqr{R0sR} and \eqr{sEACR} have been proved to hold in numerous particular examples, the fact that there exists a class of models (the SIR-PH) which unifies hundreds of previous papers has gone unnoticed until \cite{AABBGH}. \eeR

 We emphasize again that the equivalence of the two approaches for computing  $R_0$, and the formulas \eqr{R0sR}, \eqr{sEACR}, \eqr{R0fr} was proved  in \cite{AABBGH} only for SIR-PH-FA  epidemic models with one type of susceptible, and with $F$ and $K$ of rank one.

 \beR The   formula   \eqr{sEACR} for \SPF\ points to the fact that the endemic susceptibles
 are independent of the total population (the stoichiometric class), a property called \ACR\ (ACR)  in CRNT;
there, the susceptibles would  be called the \ACR\ (ACR) species.

Note that the uniqueness of the  endemic point and the formula \eqr{sEACR} above holds also  for a more refined  ``Intermediate Approximation'' of SIR-PH models introduced in (\cite{AABBGH}, Prop. 4), suggesting the following questions.
 \eeR

\beO

Are there formulas analogous to \eqr{R0sR}, \eqr{sEACR}, \eqr{R0fr}, in the case of higher rank NGMs, or several  type of susceptibles?
\eeO

\subsection{Explicit Computations for Two Examples of SIR-PH-FA Models:  SAIRS.nb, SLIAR.nb}
\beXa The SAIR/S$I^2$R model of Section \ref{s:SEIR}
 is a SIR-PH with  parameters $\va =\bep 1&0\eep, A=\bep -\g_a  &a_i\\ 0&-\g_i   \eep, \bff a= (-A)\bff 1
 =\bep a_r \\\g_i   \eep$ and
 \begin{equation}\label{sairsp}
  \bb=\bep \beta_ a \\ \beta_ i\eep, \; \mbox{so}\; B=\bep \beta_ a & 0\\ \beta_ i & 0 \eep , \bff \de=\bep 0\\\de\eep, V=\bep    \g_a   +\Lambda  &  -a_i \\
 0&  \g_i +\Lambda +\de \eep.\end{equation}

The Laplace transform of the age of infection kernel  is
 \begin{align} \label{aTr} \Hat{a}(s)&= \va (\s I+V)^{-1} \bb
 =
 \beta_ a\frac{1 }{ \left(\Lambda+\g_a  +s\right)}+
 \beta_ i \frac{a_i}{(\Lambda+\g_i  +\delta+\s ) \left(\Lambda+\g_a  +\s\right)},
 \end{align}
 and Arino et al.'s formula becomes
$
\mR= \int_0^{\I} a(\tau) d \tau= \Hat{a}(0)=\frac{\beta_ a (\Lambda+\g_i  +\delta )+a_i \beta_ i}{(\Lambda+\g_i  +\delta ) \left(\Lambda+\g_a  \right)}.$

We conclude by noting that the system \eqr{SEIRsc} admits two fixed points, with the boundary point being given by
$$\pr{\sd= \frac{\Lambda +\gamma _r}{\Lambda +\gamma _r+\gamma _s},a= 0,i= 0,\rd= \frac{\gamma _s}{\Lambda +\gamma _r+\gamma _s}}.$$

 It may be checked that the endemic point with coordinates given by
 \be{epSe}
 \bc
 s_e=\frac{1}{\mR},\\
 a_e=\frac{1}{\mR}\frac{\Lambda(\gamma_i+\delta+\Lambda)(\gamma_r+\gamma_s+\Lambda)(\mR_0-1)}{(\Lambda+\delta)(\Lambda \gamma_a+a_i \gamma_r)+\gamma_a \gamma_i\Lambda+\Lambda(\gamma_r+\Lambda)(\gamma_i+\delta+\Lambda)},\\
 i_e= \frac{1}{\mR}\frac{a_i\Lambda^2(\gamma_i+\delta+\Lambda)(\gamma_r+\gamma_s+\Lambda)(\mR_0-1)}{(\Lambda+\delta)(\Lambda \gamma_a+a_i \gamma_r)+\gamma_a \gamma_i\Lambda+\Lambda(\gamma_r+\Lambda)(\gamma_i+\delta+\Lambda)},\\
 r_e=\frac{\pp{(\gamma_a\gamma_i \Lambda+a_r\delta\Lambda+a_r\Lambda^2)(1-1/\mR)}+\pp{\frac{\gamma_s\Lambda^2+a_i \gamma_s(\Lambda+\delta)+\gamma_s\Lambda(\delta+\gamma_i)}{\mR}}}{(\Lambda+\delta)(\Lambda \gamma_a+a_i \gamma_r)+\gamma_a \gamma_i\Lambda+\Lambda(\gamma_r+\Lambda)(\gamma_i+\delta+\Lambda)},
  \ec
  \ee
   becomes positive precisely when $R_0=\sd \mR >1$, and that it is always stable when it exists (see notebook \url{}).
\eeXa
\beXa
The SIR model
 is also an  \SPF\ with  parameters $\va =\bep 1\eep, A=\bep -\g   \eep,$
  $\bb= B= \beta, V=\bep    \g   +\La \eep.$

\eeXa

\beXa The SLIAR/SEIAR epidemic model, {where $L/E$ refer to the latently infected individuals  (i.e., those who are infected but have not yet developed any symptoms)} \cite{YangBrauer,arino2020simple,AAB24} is defined by
\vspace{-9pt}
\begin{equation}\label{SLIARG}
\Scale[0.9]{ \bc
\s'(t)=   \Lambda -\s(t)\pr{\beta_ 2  i_2(t)+\beta_3 i_3(t)+\Lambda}\\
\bep i_1'(t)&i_2'(t)& i_3'(t)\eep = \bep i_1(t)&i_2(t)& i_3(t)\eep
\pp{\s(t) \bep  0    &  0& 0 \\
 {\beta_ 2}&  0 &0\\
 \beta_ 3&0&0\eep+ \left(
\begin{array}{ccc}
 -\g_1 -\Lambda & \g_{1,2} & \g_{1,3} \\
 0 & -\g_2 -\Lambda  & \g_{2,3} \\
 0 & 0 & -\g_3 -\Lambda \\
\end{array}
\right)}
\\
{\mathsf r}'(t)= \g_{2,r}  i_2(t)+\g_3  i_3(t) -\Lambda {\mathsf r}(t)\ec.}
\end{equation}
See the Figure \ref{SLIAR2} below which illustrates its corresponding shematic graph.

This is a  \SPF\  with  parameters
\vspace{-9pt}
 \bea
 \Scale[0.9]{ \va =\bep 1&0&0\eep, A= \left(
\begin{array}{ccc}
 -\g_1  & \g_{1,2} & \g_{1,3} \\
 0 & -\g_2   & \g_{2,3}\\
 0 & 0 & -\g_3   \\
\end{array}
\right), \bff a= (-A) \bff 1
 =\bep 0\\ \g_{2,r}  \\\g_3  \eep, \bb=\bep 0\\ \beta_2 \\ \beta_ 3\eep, \;   \mbox{so}\; B=\left(
\begin{array}{ccc}
 0 & 0 & 0 \\
 \beta _2 & 0 & 0 \\
 \beta _3 & 0 & 0 \\
\end{array}
\right).}
\eea
The Laplace transform of the age of infection kernel  is
\vspace{-9pt}
 \bea
\Scale[1.1]{ \Hat{a}(s)=\beta _2 \frac{ \gamma _{\text{1,2}}}{\left(\Lambda+\gamma _1+\s\right) \left(\Lambda+\gamma _2+\s\right)}+ \beta_ 3 \pr{\frac{\g_{1,3}}{\left(\Lambda+\gamma _1+\s\right) \left(\Lambda+\gamma _3+\s\right)}+\frac{\gamma _{\text{1,2}} \gamma _{\text{2,3}}}{\left(\Lambda+\gamma _1+\s\right) \left(\Lambda+\gamma _2+\s\right) \left(\Lambda+\gamma _3+\s\right)}},}\eea
and the Arino et al. formula yields
$
\mR= \frac{\beta _3 \gamma _{\text{1,2}} \gamma _{\text{2,3}}+b \beta _2 \gamma _{\text{1,2}}+\beta _2 \gamma _3 \gamma _{\text{1,2}}+b \beta _3 \gamma _{\text{1,3}}+\beta _3 \gamma _2 \gamma _{\text{1,3}}}{\left(b+\gamma _1\right) \left(b+\gamma _2\right) \left(b+\gamma _3\right)}.$

\begin{figure}[H]
\includegraphics[width=13.5cm]{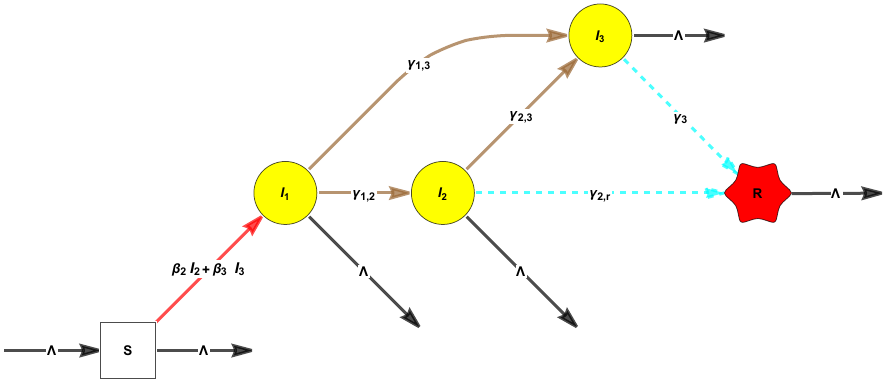}
\caption{{Chart} 
 flow/species graph of the SLIAR model \eqr{SLIARG}.
The red edge corresponds
to the entrance of susceptibles into the entrance disease class L, the cyan
dashed lines correspond to the rates of recovery, and the brown and black edges are the rates of the transition matrix V towards the interior and exterior, \resp.
}
\label{SLIAR2}
\end{figure}
\eeXa

The ODE, using slightly different notations from those in the figure, are
\vspace{-9pt}
\be{SLIAR}
\Scale[0.9]{ \bc
\s'(t)=   \Lambda -\s(t)\pr{\beta_ a  a(t)+\beta_i i(t)+\Lambda}\\
\bep {\mathsf l}'(t)&{\mathsf i}'(t)& {\mathsf a}'(t)\eep = \bep {\mathsf l}(t)&{\mathsf a}(t)&{\mathsf i}(t)\eep
\pp{\s(t) \bep  0    &  0& 0 \\
 {\beta_ a}&  0 &0\\
 \beta_ i&0&0\eep+ \left(
\begin{array}{ccc}
 -\g_l -\Lambda & l_i & l_a \\
 0 & -\g_i -\Lambda  & i_a \\
 0 & 0 & -\g_a -\Lambda \\
\end{array}
\right)}
\\
{\mathsf r}'(t)= i_r  {\mathsf i}(t)+\g_a  {\mathsf a}(t) -\Lambda {\mathsf r}(t)\ec,}
\ee
{where $\gamma_i=i_a+i_r$ and $\gamma_l=l_i+l_a.$}


The correspondence between the two sets of notation is
\bea
\bc
\gamma_1\to \gamma_l, \gamma_{1,2}\to l_i, \gamma_{1,3}\to l_a,\\
\gamma_2\to \gamma_i, \gamma_{2,3}\to i_a, \gamma_{2,r}\to i_r, \gamma_3\to \gamma_a,\\
\beta_2\to\beta_i, \beta_3\to \beta_a,\\
\pr{i_1,i_2,i_3}\to \pr{{\mathsf l},{\mathsf i},{\mathsf a}},
\ec
\eea

Proposition \ref{p:ren} yields that
\be{RSl}
\mR= \Hat{a}(0)=\frac{\beta_i l_i(\Lambda+\gamma_a)+\beta_al_a(\Lambda+\gamma_i)+\beta_a l_i i_a}{(\Lambda+\gamma_l)(\Lambda+\gamma_i)(\Lambda+\gamma_a)}\ee

We conclude by noting that the system \eqr{SLIAR} admits two fixed points, with the boundary point being given by
$$\pr{\sd= 1,a= 0,i= 0,\rd= 0}.$$

 It may be checked that the endemic point with coordinates given by
 \bea
 \bc
 s_e=\frac{1}{\mR },\\
 l_e=\frac{\Lambda}{\Lambda+\gamma_l}(1-\frac{1}{\mR }),\\
 i_e=\frac{\Lambda l_i}{(\gamma_i+\Lambda)(\gamma_l+\Lambda)}(1-\frac{1}{\mR }),\\
 a_e=\frac{\Lambda((\Lambda+i_r) l_a+\gamma_l i_a}{\beta_i l_i(\Lambda+\gamma_a)+\beta_al_a(\Lambda+\gamma_i)+\beta_a l_i i_a}(\mR -1),\\
 r_e=\frac{i_r l_i \Lambda+\gamma_a\gamma_l\gamma_i+\gamma_a l_a \Lambda}{\beta_i l_i(\Lambda+\gamma_a)+\beta_al_a(\Lambda+\gamma_i)+\beta_a l_i i_a}(\mR -1)
  \ec
  \eea
   becomes positive precisely when $R_0=\sd \mR >1$, and that it is always stable when it exists (see notebook \url{}).

The reactions of the corresponding CRN, assuming ``mass-action form'', are

$
 {0 \xrightarrow{\Lambda} S},
 {L + S \xrightarrow{a s \beta _a}  2 L},
 {I + S \xrightarrow{ i s \beta _i}  I + L},
 {L \xrightarrow{l l_a} A},
 {A \xrightarrow{a \g_a} R},
 {I \xrightarrow{a i i_a } A},\\
 {L \xrightarrow{l l_i} I},
 {I \xrightarrow{i i_r} R},
 {S \xrightarrow{\Lambda  s} 0},
 {L\xrightarrow{\Lambda  l} 0},
 {I \xrightarrow{\Lambda i} 0},
 {A \xrightarrow{\Lambda a} 0},
 {R \xrightarrow{\Lambda  r} 0},
$

see the Figure \ref{SLIARgr} below.

\begin{figure}[H]
\includegraphics[width=13cm]{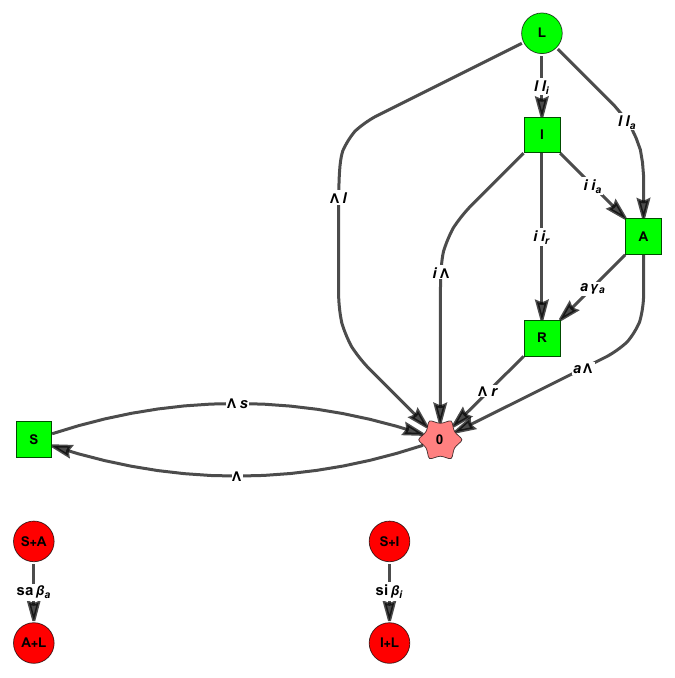}
\caption{{The} 
 {Feinberg--Horn--Jackson graph} of the ``SLIAR network'' with $13$ reactions. The deficiency is $10-3-5=2$, and weak reversibility does not hold.}
\label{SLIARgr}
\end{figure}

\subsection{SIR-PH-FA Meets ACR: EnvZ.nb}\lab{s:SACR}

We provide here an ACR example, for which the stability condition for the fixed boundary point is of the  form \eqr{R0fr}.

Consider a simplified EnvZ-OmpR  system \cite{AAHJ}:
\begin{equation}
\label{envzompr}
\begin{tikzcd}
X  \arrow[r,"k_1"] & X_t  \arrow[r,"k_2"] & X_p \\[-0.1in]
X_p + Y \arrow[r,"k_3"] & X  + Y_p \\[-0.1in]
X_t  + Y_p \arrow[r,"k_4"] & X_t  + Y
\end{tikzcd}
\end{equation}
with \sm
$\left(
\begin{array}{cccc}
 -1 & 0 & 1 & 0 \\
 1 & -1 & 0 & 0 \\
 0 & 1 & -1 & 0 \\
 0 & 0 & -1 & 1 \\
 0 & 0 & 1 & -1 \\
\end{array}
\right),$ rank $3$,  deficiency $7-3-3=1$, and a corresponding system of ordinary differential equations:

\begin{equation}
\label{de}
\left\{ \begin{aligned}
X' & = -k_1 X  + k_3 X_p Y \\
X_t ' & = k_1 X  - k_2 X_t \\
X_p' & = k_2 X_t  - k_3 X_p Y \\
Y' & = -k_3 X_p Y + k_4 X_t  Y_p \\
Y_p' & = k_3 X_p Y - k_4 X_t  Y_p
\end{aligned} \right.
\end{equation}
(to simplify the notation, we  have abbreviated, as is customary,  EnvZ by X, and OmpR by Y). There are two conservation laws $x_{tot} = X  + X_t  + X_p, y_{tot} = Y + Y_p$, (total EnvZ and total OmpR, \resp).

  Mathematica finds quickly the ``DFE'' \fp\ with resident species  $\{ X_p, Y_p\}$:
  $$X  =  X_t  =0=Y,  X_p  =x_{tot},  Y_p =y_{tot},$$
 and possibly a second interior \fp, parametrized by $(x  ^*, x_p  ^*)$: \[x_t  ^*=\frac{k_1 x  ^*}{k_2}, y  ^*= \frac{k_1 x  ^*}{k_3 x_p  ^*}, y_p ^*=\frac{k_2}{k_4} \leq y_{tot},\]
the last component $y_p ^*$ being the ACR species
(this example  satisfies the sufficient \cite{ShiFei}
ACR conditions since it  has deficiency one and two non-terminal complexes
$X_t, X_t+Y_p$---see \eqr{envzompr}---which
  differ in the single species $Y_p$, which must therefore be  the ACR species).

It turns out that the last obvious inequality implies  both the uniqueness of the interior point within each invariant set, parametrized  by $x_{tot}, y_{tot}$, its stability, and the instability of the DFE when the inequality is strict. Rewriting this inequality in the form
\be{R0A} 1 < \fr{y_{tot}}{y_p ^*}:=R_0,\ee
which is the complement of the inequality $R_0\leq 1$ \eqr{R0fr}, we render obvious its CRN interpretation that for the interior point to be stable, it must first satisfy the total constraint where it intervenes.

\beR \lab{r:hNGM} Let us also establish the stability of the DFE,
as a full example of the NGM approach, in the non-rank one case (so, the formula of \cite{Arino} does not apply).
\BEN \im The derivative of $RHS_i$ \wrt\ the infection variables is
\begin{gather}\label{jI} \nonumber J_i=\left(
\begin{array}{ccc}
 -k_1 & 0 & k_3 X_p   \\
 k_1 & -k_2 & 0 \\
 0 & k_4 Y_p   & -k_3 X_p   \\
\end{array}
\right)=F-V=\left(
\begin{array}{ccc}
 0 & 0 & k_3 X_p   \\
 0 & 0 & 0 \\
 0 & k_4 Y_p   & 0 \\
\end{array}
\right)-\\
\left(
\begin{array}{ccc}
 k_1 & 0 & 0 \\
 -k_1 & k_2 & 0 \\
 0 & 0 &
 k_3 X_p   \\
\end{array}
\right) \Lra
\end{gather}

\be{K} K = FV^{-1} = \left( \begin{array}{ccc}
0 & 0 & 1 \\
0 & 0 & 0 \\
\displaystyle{\frac{k_4 y_{tot}}{k_2}} & \displaystyle{\frac{k_4 y_{tot}}{k_2}} & 0
\end{array}
\right)
\ee
where the resident species $X_p, Y_p$ are evaluated
at the ``disease-free state'' (with $X  = X_t = Y = 0$),
i.e., $X_p = x_{tot}, Y_p = y_{tot}$ (this is implemented in our script NGM).
\im The next step consists in attempting to
factor  the \ch\  of $K$, and in removing ``stable factors'' (i.e., factors which may not have \eig s with a positive real part). It seems, at least for the ACR models \satg\ the \cite{ShiFei} conditions discussed in this section,   that the \ch\ of $K$ factors more often than that of $J_i$.

 In our current example, the 3rd degree \ch\ of $K$ is $ch(u)=u (k_2 u^2-k_4 Y_p)$. After removing the factor $u$,    the remaining second-degree polynomial has symmetric roots, the positive one being
$$R_0=\sqrt{\frac{k_4 y_{tot}}{k_2}}=\sqrt{\frac{ y_{tot}}{y_p  ^*}}.$$

Finally, the condition  $R_0 > 1$ under which the stability of the boundary  disease-free steady state is lost is precisely
\be{ACRi}y_{tot} > \frac{k_2}{k_4}=y_p^*.\ee
\EEN
\eeR

\beR When the  \ch\ has a higher degree, we may also apply,  alternatively, the RH criteria to the shifted polynomial $p(x)=ch(x+1)$. In our current example,
the shifted second-degree polynomial is $k_2 x^2+2 k_2 x+k_2-k_4 Y_p$, and RH recovers
the condition $k_2-k_4 Y_p\ge 0$. The computational advantage of this approach is that the single condition $R_0 \leq 1,$  with $R_0$, which may be the root of a high-order polynomial, may be replaced by several RH conditions.
\eeR

 \beO Are varying total population models and their intermediate approximations of interest in CRN? Are there particular classes of CRN models where formulas similar to \eqr{sEACR} and  \eqr{R0fr} hold (which include for example the ACR examples provided in \cite{AAHJ})?
 \eeO

\beR As an aside,    nowadays, \ME\ literature suffers   from "irreproducibility'', i.e., the lack of electronic notebooks to support complicated computations, an absence which changes the  simple task of pressing enter for checking into days of work.

To remedy this, we have provided our Mathematica package and notebooks for the above{examples at} 
 \url{https://github.com/adhalanay/epidemiology_crns } (accessed on 18 September 2024). Since they may still contain mistakes or unclarities occasionally, we ask the reader to contact us for any remarks. \eeR

\section{Further Fundamentals of CRNT: The Euclidean FHJ Graph, Network Translation GMAKs, and ACR Results; The (4,5,3,2) MAK (\cite{tonello2018network}, (1)): Tonello.nb}\lab{s:Ton}

We end our essay by mentioning one idea which seems  important for the non-reversible models of ME. The reader may have noticed that one idea
of CRNT, the embedding of the \FHJ\ into Euclidean space \cite{CF05}, raises the question  of  how moving the graph vertices, and allowing them maybe to collide, will affect the corresponding ODE. This implicit question has been now addressed in the theory of GMAKs obtained by network translation, initiated by Johnston  \cite{John,tonello2018network,JMP,hernandez2022framework}. This approach requires understanding that translation has a double effect: (1) it changes the directions in the SM, and (2) it also changes the matrix of exponents $Y_{\al}$---see Example \ref{e:linSIRS}. An interesting effect is obtained when we effectuate only the first change but leaving $Y_{\al}$ and $diag({\bff k})$ unchanged.
It can easily be seen that the  SM and ODE for the resulting GMAK remain then unchanged, provided all complexes in a linkage class are translated by the same vector. This way, one may achieve weak reversibility and ``kinetic deficiency'' $0$ for the new SM, as happened in Example \ref{e:linSIRS}.

 Network translation of a MAK  to a ZD-WR GMAK  helps understanding, among other things,  rational parametrizations related to the  matrix tree theorem (see, for example, \cite{CoxSl}) which appear in certain MAKs. This is well illustrated by  the MAK  (\cite{tonello2018network}, (1)) \begin{equation}\label{eq:proper_acr_intro}
      \begin{tikzcd}[row sep=small, column sep=small]
          & A + C & C \arrow[rd, "r_3"] & & \\
         A + B \arrow[ru, "r_1"] \arrow[rd, "r_2"'] & & & A \arrow[r, "r_5"] & B. \\
        & A + D & D \arrow[ru, "r_4"'] & &
      \end{tikzcd}
\end{equation}
with \sm\ $\left(
\begin{array}{ccccc}
 0 & 0 & 1 & 1 & -1 \\
 -1 & -1 & 0 & 0 & 1 \\
 1 & 0 & -1 & 0 & 0 \\
 0 & 1 & 0 & -1 & 0 \\
\end{array}
\right)$.

Ordering  the species (and the complexes) according to their appearance order yields the ODE: \beq &&\lab{ODET}\frac{d}{dt} \bep A\\
 B\\
 C\\
   D\eep=
   \bep-A k_5+C k_3+D k_4\\-A \left(B k_1+B k_2-k_5\right)\\
   A B k_1-C k_3\\A B k_2-D k_4\eep\eeq

\BEN \im
It can be checked that $rank(\Gamma)=3$, the deficiency is
$\delta=n_C-rank(\Gamma)-\ell=7-3-2=2,$
 the CRN has the conservation
$
A+B+C+D
$, and that the cone of positive fluxes has dimension 2.

 \im \Ma\ finds two fixed points
\be{Ton}\bc A= 0,C= 0,D= 0,B=n_{tot}\}\\\left\{b^*= \frac{k_5}{k_1+k_2},c^*= a^*\frac{ k_1 k_5}{\left(k_1+k_2\right) k_3},d^*= a^*\frac{ k_2 k_5}{\left(k_1+k_2\right) k_4}\right.\ec,\ee
which suggests that the \MTT\ might be at work.
The second solution under the extra constraint $A+B+C+D=n_{tot}$,
\[\bc A=\left(\left(k_1+k_2\right) n_{tot}-k_5\right)\frac{k_3 k_4 }{k_1 k_4 \left(k_3+k_5\right)+k_2 k_3 \left(k_4+k_5\right)},B= \frac{k_5}{k_1+k_2}\\C=\left(\left(k_1+k_2\right) n_{tot}-k_5\right)\frac{k_1 k_4 k_5 }{\left(k_1+k_2\right) \left(k_1 k_4 \left(k_3+k_5\right)+k_2 k_3 \left(k_4+k_5\right)\right)}\\D=\left(\left(k_1+k_2\right) n_{tot}-k_5\right)\frac{k_2 k_3 k_5 }{\left(k_1+k_2\right) \left(k_1 k_4 \left(k_3+k_5\right)+k_2 k_3 \left(k_4+k_5\right)\right)}\ec\] is
positive {if and only if} 
$n_{tot} >b^*$.
The NGM method yields $R_0=\fr{n_{tot}}{b*}$, and the Jacobian at the EE reveals that it is always stable when it exists, i.e., when $R_0>1$.
\EEN

To better understand this model, the authors propose studying in parallel the GMAK
\begin{equation}\label{eq:proper_acr_intro_translated}
      \begin{tikzcd}[row sep=small, column sep=small]
        & & \mbox{\ovalbox{$\begin{array}{c} A+C \\ (C) \end{array}$}} \arrow[rd, "\tilde{r}_3"] & & & \\
& \mbox{\ovalbox{$\begin{array}{c} A+B \\ (A+B) \end{array}$}} \arrow[ru, "\tilde{r}_1"] \arrow[rd, "\tilde{r}_2"'] & & \mbox{\ovalbox{$\begin{array}{c} 2A \\ (A) \end{array}$}}. \arrow[ll, "\tilde{r}_5"] \\
        & & \mbox{\ovalbox{$\begin{array}{c} A+D \\ (D) \end{array}$}} \arrow[ru, "\tilde{r}_4"'] & & &
      \end{tikzcd}
    \end{equation}

    In \eqref{eq:proper_acr_intro_translated}, the lower (kinetic) complexes, displayed in parenthesis,  correspond to the original source complexes \cite{John}, and the upper (stoichiometric) complexes correspond to the source complexes adjusted by  ``translation complexes'' (in this case,   $A$ is added to $r_3$, $r_4$, and $r_5$), with the net result of gluing some of the product complexes to source ones, creating in this way a WR GMAK.

    Now the translated CRN, which represents the same dynamics, has ZD, and hence the \fp\ formulas may be interpreted via the \MTT, applied to the graph \eqref{eq:proper_acr_intro_translated} (\cite{tonello2018network}, Thm 2-3).

We end this section by showing how to  use the Matlab TowardZ algorithm \cite{HongSon,Hong}, which uses
network translation,  a{vailable  at}  \url{ https://github.com/Mathbiomed/TOWARDZ}, (accessed on 5 December 2022), for determining some WR-ZD network translations for the  MAK \eqr{eq:proper_acr_intro}, including that given in \eqr{eq:proper_acr_intro_translated}.
\BEN \im  We must type first the sources and products for the MAK of (\cite{tonello2018network}, (1)) in the file TOWARDZ\_for\_given.m, which is provided in  TowardZ for testing the illustrated examples from the paper (or others).

\im The  results to examine first are the solution, a 7 {$\times$} 
 2 cell array object, and Index, a 7 {$\times$} 5 cell array. {They indicate the presence of seven distinct WR-ZD GMAK realizations. The source and product for the GMAK \eqr{eq:proper_acr_intro_translated} may be found by typing  Solution\{7,:\} .} This  corresponds to the solution provided in (\cite{tonello2018network} (1)), after a permutation of the reactions provided by Index \{7,:\} (which reveals that  the source and product complex matrices displayed have the column order r4, r3, r2, r1, and r5).
\EEN

\subsection{Finding  Weakly Reversible and
Zero Deficiency (WR-ZD/$WR_0$) Representations of an ODE,  Using TowardZ \cite{Hong}}\label{s:tow}

Since the existence of a WR-ZD representation guarantees that an ODE has remarkable properties, it is quite interesting to detect the existence of such representations. Let us first mention one algorithm for that, which   offers MAK representations,  provided  in the recent paper \cite{CraJinYu}. We have not yet succeeded in finding ME examples where
it applies, so we state this as  an open problem.
    \beO
    Are there ME models which admit equivalent MAK WR-ZD representations? \eeO

  \beR   This question is related to the papers \cite{PuenteJohn,BurtonJohn,hoessly2023complex,CraJinYu}, where the question of which models admit WR-ZD representations is  studied algorithmically but may not be fully resolved theoretically.
  \eeR

 On the other hand,  TowardZ  does produce  GMAK representations for several of the ME models we tried.

We provide now the 4 TowardZ WR-ZD translations of the deficiency one SIRS model with eight reactions \eqr{SIRS}:
\vspace{-9pt}
\[\frac{d{x}(t)}{dt}=\left(
\begin{array}{cccccccc}
 1 & -1 & 0 & 1 & -1 & -1 & 0 & 0 \\
 0 & 1 & -1 & 0 & 0 & 0 & -1 & 0 \\
 0 & 0 & 1 & -1 & 1 & 0 & 0 & -1 \\
\end{array}
\right) \bep\lambda \\ \beta  i s\\ \gamma  i\\ r \gamma _r\\ s \gamma _s\\ \mu  s\\ i \mu _i\\ \mu  r\eep=\left(
\begin{array}{c}
- \beta  s i+\lambda +r \gamma _r-s \gamma _s-\mu  s \\
 -\gamma  i-i \mu _i+\beta  i s \\
 \gamma  i-r \gamma _r-\mu  r+s \gamma _s \\
\end{array}
\right)\]

\beR \lab{r:tow} Besides the monomolecular SIRS from example \ref{e:linSIRS}, TowardZ finds three other ``WR-ZD cousins'' SIRS:
\BEN \im \{{"r'' $\to$ "r'' + "s'', "r'' + "s'' $\to$ "i'' + "r'', "i'' + "r'' $\to$ "r'',
 "r'' + "s'' $\to$ 2 "r'', 2 "r'' $\to$ "r'' + "s'', "r'' + "s'' $\to$ "r'',
 "i'' + "r'' $\to$ 2 "r'', 2 "r'' $\to$ "r"}\}

 \im  \{"i'' $\to$ "i'' + "s'', "i'' + "s'' $\to$ 2 "i'', 2 "i'' $\to$ "i'',
 "i'' + "s'' $\to$ "i'' + "r'', "i'' + "r'' $\to$ "i'' + "s'', "i'' + "s'' $\to$ "i'',
 2 "i'' $\to$ "i'' + "r'', "i'' + "r'' $\to$ "i"\}

 \im \{"s'' $\to$ 2 "s'', 2 "s'' $\to$ "i'' + "s'', "i'' + "s'' $\to$ "s'',
 2 "s'' $\to$ "r'' + "s'', "r'' + "s'' $\to$ 2 "s'', 2 "s'' $\to$ "s'',
 "i'' + "s'' $\to$ "r'' + "s'', "r'' + "s'' $\to$ "s"\}
 \EEN

Now SIRS is in itself an easy model, but it is quite possible that such cousins might turn out useful for studying more complex models.
\eeR

\section{Discussion: Can  (Generalized) Chemical Reaction Network  Methods Help in Solving Mathematical Epidemiology Problems, and Vice Versa?
} \lab{s:SIRCRN}
The CRN formalism provides a universal language for studying essentially non-negative ODEs, which is not known well enough outside CRNT;  as a result, classic results like, for example, that in \cite{hun}, and the use of reaction variables \cite{Clarke} are  rediscovered again and again, while newer results
  like \cite{banaji2018inheritance} are of course ignored.

  On the other hand, ME models are  quite difficult to analyze under realistic assumptions, as witnessed by the huge proportion of ME papers that work under approximations like the negligibility of deaths, permanent immunity,  single-strain viruses, and either linear birth rates or constant  immigration  rates (but not a combination of both).
  Let us stress that the simplest ME model, the three species SIRS, has only been analyzed under complete realistic assumptions only recently, in \cite{nill23}, and that despite brilliant papers on SEIRS like, for example, \cite{LiMul,LiGraef,LiMulVan,SunHsieh,LuLu}, this four-species model is far from being fully understood.
  Also, as already noted, simplifying properties in CRN, like weak reversibility and ZD, seem to be never met in  ME models (Open Problem 1).

 Therefore, the answer to  the question in the title is not easy. What has been established for now  is that the intersection of ``classical chemical models''  and ``ME-type models with conditionally stable boundary fixed points'', while not large, is not empty since it includes the CRN ACR models, for which  the NGM method for establishing DFE stability turns out to be quite convenient, and for which relations first encountered for SIR-PH-FA models like \eqr{R0sR} hold.

\Fr  as recalled already in Section \ref{s:syn},  two CRN methods found spectacular ME applications in \cite{AAV}, namely, the following: \BEN \im
 The inheritance of Hopf bifurcations was applied to  SIRnS models {(i.e., cyclic epidemic models with one susceptible class, one infection class, and $n$ recovered classes)}.
  \im The existence of bifurcations for models with rich rates  has   turned  out useful  for showing that all models which include as submodels a Capasso-type SIR  (i.e., SIR models with admissible symbolic incidence function) admit Hopf bifurcations \cite{AAV}.\EEN

   Our personal answer to the question in the title of this section is that the CRN formulations and notation style seem more suitable to us for capturing ``network dynamics'' than the traditional  representation.  For that, we  use them nowadays
  daily in investigating any essentially non-negative ODE problem, by calling our CRN-style Mathematica package EpidCRN, {offered at} \url{ https://github.com/adhalanay/epidemiology_crns}, (accessed on 18 September 2024), (already used in~\cite{AAHJ}), which, in turn, uses the package ReactionKinetics.

But a complete answer to the question in the title is impossible since CRN/BIN theory is  quickly evolving nowadays. For some developments which might turn out useful also in ME, let us suggest the theory of k-contractions~\cite{WeissM21,kcon,AAS23str,RamiM,ARS24}, of  structural stability, and BDC decompositions~\cite{BFG14,BG15,GiuliaFranco,Giulia,Colaneri,BG18}, and that of robust Lyapunov functions~
\cite{BG14,AAunc,AASrob,AAnew,Ali,BGS18,Ali23gr}.

In the opposite direction,
 the renewal kernel method---see for example \cite{Diek,AAHJ}, and Section \ref{s:renker} above---might turn out to be useful  for studying CRN  models.

 Let us end this paper by repeating the question raised in Remark \ref{r:tow}:  can the ``WR-ZD cousins'' of epidemic models shed some light on their ME relatives?

 \section{{A Brief} Tour  of Some Relevant Facts from the Theory of ``Continuous Time Pure Jump Markov Chains'' (CTMC)}
 \lab{s:CTMC}

An ODE  RHS like \eqr{SYRPH} may also be used  to define a stochastic CTMC (continuous-time Markov chain) process on a mesh $h \N^3,h >0$, which jumps in the direction of the columns of $\Gamma$, at  rates $\bf r(x)$.
 Note that for the Mathematica kernel/ChatGPT/Alexa,
the only difference between the specification of an ODE and of its associated pure jump Markov model  is an extra ``mesh parameter'' $h$.

Since only numeric results are available typically for multidimensional CTMC processes, various semi-analytic
approximations have been proposed as well, notably by fixing the non-infectious compartments, since intuitively they evolve on  a smaller time-scale.

Citing Griffiths \cite{griffiths1973multivariate},
``It has been noted by Bartlett (1955), p. 129, that for an epidemic in a large
population, the number of susceptibles may, at least in the early stages of an
outbreak, be regarded as approximately constant at its initial value and that this
approximation will continue to hold throughout the course of an epidemic,
provided that the final epidemic size is small relative to the total susceptible
population. Thus the general (SIR) epidemic process may be approximated by a one-dimensional
birth-and-death process."

The resulting approximation may either  converge to $0$ or  be non-positive--recurrent (due the infinite state space) as detailed in the next section.

\subsection{An Approximate ''Stochastic SIR Infection Process" on the infected compartment}
An approximate ``stochastic SIR infection process'' obtained by fixing $s,r$ is associated to the equation $$i'=\beta s i-\mu i=(\beta s -\mu) i,$$ of the disease \com\ $i$, where $s$ is assumed fixed. Recall that this projection on the disease \com s is also the main brick of the next generation matrix (NGM) method for computing the stability of the DFE \cite{Diek,Van,Van08}.

  The resulting birth and death process   is a Markov process $X_t \in \N$
  with {\bf linear rates},   generating operator on the set of  functions $f:h\N \to h\N$ defined by
\bea
\mG f(i)= \beta s i (f(i+h)-f(i)) -\mu   i( f(i-h)-f(i)):=A f(i).\eea
This process either converges to $0$ or  is non-positive--recurrent, depending on whether
$R_0:=\fr {\beta s}{\mu } $ is strictly smaller than $1$, or not. The probabilities of "extinction/absorption into $0$'',
when starting  the process with  $j$ infected are
\be{q} q^j, \qu q=\bc 1& R_0 \leq 1\\ \fr {\mu } {\beta s}=\fr 1{R_0} & R_0 > 1 \ec\ee
see, for example, the textbook \cite{dawson2017introductory}.

A general quadratic (matrix) formula for computing the extinction probabilities for any epidemic model inspired by Bacaer and spelled out in \cite{AABJ} is also implemented in our package EpidCRN.

\subsection{A Review of the Times of Absorption of Finite-State Continuous-Time Markov Chains}\la{s:PH}

There exists a second approximate epidemic process of interest, for which there does not exist much work in the literature (probably due to its trivial nature).  Recall that we fixed all the non-infection \com s so that now all the transition rates except for birth are linear. Assume now we are following the evolution of a single infected individual between the infection \com s, which changes the rates from linear to constant, given by the elements of the matrix $-V$. We have now to deal with a Markovian evolution among a set of finite states, i.e., a finite state CTMC, and the issue is to study its time until absorption (which corresponds to the infected individual becoming recovered.

 A phase-type   distribution is the distribution of the time until the absorption of a finite-state continuous-time Markov chain (CTMC) into an absorbing state. More formally, we have the following.

\beD[\cite{Ramaswami99}]
Let $X_t$ denote a finite-state continuous-time Markov process, with one   absorbing state.
  Denote by $\bep \bold{a}& A\\ 0 & \bold{0}\eep$  the generating matrix of $X_t$, where $\bold{a}:=-A {\bf 1}$, and  $A$ is an $n\times n$ matrix, describing transitions between transient states, and let ${\vec \alpha}=\pr{\alpha_1,\alpha_2,\dots ,\alpha_n}$ denote a 1-dimensional subprobability row vector of size $n$, representing the initial distribution of $X_t$.

 The ``absorption time''  until $X_t$ enters its absorbing state, to be denoted by $T_{abs}$, is also called a phase-type random variable.

 \eeD

The matrix-exponential $e^{t A}$ yields, by definition, the transition probabilities at time $t$ between the transient states. The following hold:
\BEN \im
The survival function of $T_{abs}$ is given by
\be{PHs} P[T_{abs}>t]={\vec \alpha} e^{t A} {\bf 1},\ee
where $ {\bf 1}$ is a column vector of $1$'s.
\im The density of $T_{abs}$ is
\be{PHd}f_{T_{abs}}(t)={\vec \alpha} e^{t A} {\bf a},\ee
where $ {\bf a}=(-A) {\bf 1}$ is a column vector representing the direct absorption rates (via one transition).
\im The ``expected dwell times'' in each state are given by the components
of the matrix $(-A)^{-1}$:
\be{edt}E_i[1_{\{X_t=j\}}] =(-A)^{-1}_{i,j}\ee

\EEN
\beR The expression above shows that the elements of the NGM matrix \eqr{K} may be interpreted as the total force of infectivity an infectious individual will exert before recovery. \eeR

The expression phase-type   distribution (or law) refers to the ensemble of probabilities $P[T_{abs} \in [a,b], a, b \in \R_+$, and sometimes to just one of the two equations, \eqr{PHs} and \eqr{PHd}, given above.

 Phase-type   distributions have become very popular in applied probability,
 due to their computational tractability. Essentially, all known results on exponential distributions hold for these more complex laws (simply by replacing scalars by matrices or vectors). Since this class is dense in the class of continuous laws on $R_+$, statisticians may always fit PH distributions to their data, while simultaneously taking advantage
 of the exceptional  tractability of the exponential law.

The ancestor of the phase-type modelization is the method of stages,  introduced by A. K. Erlang. In ME, phase-type modelization was popularized by Hurtado \cite{Hurtado19,hurtado2021building}---see also \cite{AAB,AAK} for some explanation of why probabilities intervene in  ODE  ME models.

\vspace{12pt}
\subsection{The CRN representation yields a Hamiltonian representation of the expected time spent along ``escape paths" of  stochastic CTMC   processes}\lab{s:esc}

{

The importance of the CRN formalism becomes apparent  when introducing the
 Hamiltonian function associated to a CRN (see for example \cite{Snarski,Gagrani,clancy2024extinction}):
 }
{
\begin{align}
    H_\text{CRN}(\theta,x)
    &=
    \sum_{r \in \mR}\big( e^{(y_{\beta}-y_{\alpha})\cdot \theta} - 1\big) k_{y_{\alpha}\to y_{\beta}}x^{y_{\alpha}} \label{eq:Ham_CRN_rxn}
\end{align}
}
{
(note that this depends only on the stoichiometric columns $y_{\beta}-y_{\alpha}$, and that when $\theta \to 0, H_\text{CRN}(\theta,x)/\theta \to RHS$ ).
The Hamiltonian allows computing the trajectories of the least improbable paths by which a CTMC escapes from a stable point (say endemic) towards an unstable one (say the DFE), but the rarity of such events  makes them uninteresting in ME.}


A Hamiltonian has an associated  deterministic motion in $\left( x , \theta\right)$ space, defined by the equations:
\begin{eqnarray} \label{Hamiltonian_system}
{dx \over dt} = {\partial H \over \partial \theta} , \hspace{2cm} {d \theta\over dt} = - \,  {\partial H \over \partial x} .
\end{eqnarray}
{
The expected travel time to an unstable deterministic fixed point,  starting from a stable deterministic equilibrium point,  satisfies (see~\cite{EK04}) \be{action}\fr{\ln \left( \tau_{\lfloor N x^* \rfloor} \right)} N \to A=\int_{-\infty}^\infty \theta\, {dx \over dt} \, dt \text{ as }N \to \infty,\ee where $A$ is known as  ``action functional".
The value of $A$ is found by integrating along the trajectory of~(\ref{Hamiltonian_system}) that goes from initial state~$(x,\theta) = \left( x^* , 0 \right)$ at time $t = -\infty$, to a final state at time $t=+\infty$ which is a fixed point of the Hamiltonian ~$(x,\theta) = \left( 0 , 0\right)$.

The methods outlined above are quite general, and can potentially be applied numerically to any ME model.

As an aside, for the SI process  the value of $A$ may be  found  explicitly \cite{DSS05,Clancy}, but for higher dimensions only numeric integration seems to work.
Whether or not symbolic treatment is possible in other particular cases is an open problem.

\bibliographystyle{amsalpha}
\bibliography{Pare40}

\end{document}